\documentclass[12pt]{amsart}
\usepackage[dvips]{epsfig}
\usepackage{graphics}
\usepackage{latexsym}
\usepackage{verbatim}
\usepackage{amsmath}
\usepackage{amsthm}
\usepackage{amssymb}
\newtheorem{theorem}{Theorem}[section]
\newtheorem{lemma}[theorem]{Lemma}

\newtheorem{proposition}[theorem]{Proposition}
\newtheorem{definition}[theorem]{Definition}

\def\Remark{\medskip\noindent{\bf Remark: }}
\def\Remarks{\medskip\noindent{\bf Remarks: }}

\newcommand{\ens}[1]{\mathbb{#1}}
\newcommand{\ron}[1]{\mathcal{#1}}
\newcommand{\N}{\mathbb{N}}

\newcommand{\R}{\mathbb{R}}
\newcommand{\C}{\mathbb{C}}

\def\cal{\mathcal}

\def\derpar#1#2{\frac{\partial#1}{\partial#2}}
\def\var{\varepsilon}
\def\signcm{\bigskip\bigskip\hspace{80mm}
\vbox{{\sc C. Mouhot\par\vspace{3mm}
UMPA, \'ENS Lyon\par
46 all\'ee d'Italie\par
69364 Lyon Cedex 07\par
FRANCE\par\vspace{3mm}
e-mail:} cmouhot@umpa.ens-lyon.fr }}

\begin{document}

\title[Rate of convergence to equilibrium for the Boltzmann equation]
{Rate of convergence to equilibrium for the spatially 
homogeneous Boltzmann equation with hard potentials}

\author{Cl\'ement Mouhot}

\hyphenation{bounda-ry rea-so-na-ble be-ha-vior pro-per-ties
cha-rac-te-ris-tic}

\begin{abstract} 
For the spatially homogeneous Boltzmann equation with hard potentials 
and Grad's cutoff (e.g. hard spheres), we give quantitative 
estimates of exponential convergence to equilibrium, and 
we show that the rate of exponential decay is governed 
by the spectral gap for the linearized equation, on which 
we provide a lower bound. 
Our approach is based on establishing  
spectral gap-like estimates valid near the 
equilibrium, and then connecting the latter to the quantitative nonlinear theory. 
This leads us to an explicit study of the  
linearized Boltzmann collision operator in functional 
spaces larger than the usual linearization setting.
\end{abstract}

\maketitle

\textbf{Mathematics Subject Classification (2000)}: 76P05 Rarefied gas
flows, Boltzmann equation [See also 82B40, 82C40, 82D05]; 35B40 Asymptotic behavior of solutions. 

\textbf{Keywords}: Boltzmann equation; spatially homogeneous; 
hard spheres; rate of convergence to equilibrium; explicit; spectral gap; sectorial; 
entropy method. 

\tableofcontents

\section{Introduction}
\setcounter{equation}{0}

This paper is devoted to the study of the asymptotic behavior 
of solutions to the spatially homogeneous Boltzmann equation for hard potentials with cutoff. 
On one hand it was proved by Arkeryd~\cite{Ar88} by non-constructive arguments that 
spatially homogeneous solutions (with finite mass and energy) of the Boltzmann 
equation for hard spheres 
converge towards equilibrium with exponential rate, with no 
information on the rate of convergence and the constants (in fact the proof in this paper 
required some moment assumptions, but the latter can be relaxed with the results 
about appearance and propagation of moments, as can be found in \cite{Wenn:momt:97}). 
On the other hand it was proved in~\cite{MoVi} a quantitative convergence result 
with rate $O(t^{-\infty})$ for these solutions. The goal of this paper 
is to improve and fill the gap between these results by 
 \begin{itemize}
 \item showing exponential convergence towards equilibrium by constructive arguments (with explicit rate 
 and constants);
 \item showing that the spectrum of the linearized 
 collision operator in the narrow space $L^2(M^{-1}(v)dv)$ ($M$ is the equilibrium) 
 dictates the asymptotic behavior of the solution 
 in a much more general setting, as was 
 conjectured in~\cite{CaGaTo} on the basis of the study of the Maxwell case. 
 \end{itemize}
Before we explain our results and methods in more details let us introduce the problem in a precise way. 

\subsection{The problem and its motivation}
The Boltzmann equation describes the behavior of
a dilute gas when the only interactions taken into account are 
binary collisions, by means of an evolution equation 
on the time-dependent particle distribution function in 
the phase space. 
In the case where this distribution function is
assumed to be independent of the position, we obtain 
the {\it spatially homogeneous Boltzmann equation}:
 \begin{equation}\label{eq:base}
 \derpar{f}{t}  = Q(f,f), \quad  v \in \R^N, \quad t \geq 0
 \end{equation}
in dimension $N \ge 2$. In spite of the strong restriction that this assumption 
of spatial homogeneity constitutes, it has proven an interesting 
and inspiring case for studying qualitative properties of the 
Boltzmann equation. In equation~\eqref{eq:base}, $Q$ is the quadratic
{\it Boltzmann collision operator}, defined by the bilinear form
 \begin{equation*}\label{eq:collop}
 Q(g,f) = \int _{\R^N \times \ens{S}^{N-1}} B(|v-v_*|, \cos \theta)(g'_* f' - g_* f) \, dv_* \, d\sigma.
 \end{equation*}
Here we have used the shorthands $f'=f(v')$, $g_*=g(v_*)$ and
$g'_*=g(v'_*)$, where
 \begin{equation*}\label{eq:rel:vit}
 v' = \frac{v+v_*}2 + \frac{|v-v_*|}2 \, \sigma, \qquad
 v'_* = \frac{v+v_*}2 - \frac{|v-v_*|}2 \, \sigma
 \end{equation*}
stand for the pre-collisional velocities of particles which after
collision have velocities $v$ and $v_*$. Moreover $\theta\in
[0,\pi]$ is the {\em deviation angle} between $v'-v'_*$ and $v-v_*$, and
$B$ is the Boltzmann {\em collision kernel} determined by physics 
(related to the cross-section $\Sigma(v-v_*,\sigma)$ 
by the formula $B=|v-v_*| \, \Sigma$). On physical grounds, it is
assumed that $B \geq 0$ and $B$ is a function of $|v-v_*|$ and
$\cos\theta$. 

Boltzmann's collision operator has the fundamental properties of
conserving mass, momentum and energy
  \begin{equation*}
  \int_{\R^N}Q(f,f) \, \phi(v)\,dv = 0, \quad
  \phi(v)=1,v,|v|^2 
  \end{equation*}
and satisfying Boltzmann's $H$ theorem, which can be formally written as 
  \begin{equation*} 
  {\cal D}(f):= - \frac{d}{dt} \int_{\R^N} f \log f \, dv = - \int_{\R^N} Q(f,f)\log(f) \, dv \geq 0.
  \end{equation*}
The $H$ functional $H(f) = \int f \log f$ is the opposite of the entropy of the solution.
Boltzmann's $H$ theorem implies that any equilibrium distribution
function has the form of a Maxwellian distribution
  \begin{equation*}
  M(\rho,u,T)(v)=\frac{\rho}{(2\pi T)^{N/2}}
  \exp \left( - \frac{\vert u - v \vert^2} {2T} \right), 
  \end{equation*}
where $\rho,\,u,\,T$ are the density, mean velocity
and temperature of the gas
  \begin{equation*}
  \rho = \int_{\R^N}f(v) \, dv, \quad u =
  \frac{1}{\rho}\int_{\R^N}vf(v) \, dv, \quad T = {1\over{N\rho}}
  \int_{\R^N}\vert u - v \vert^2f(v) \, dv, 
  \end{equation*}
which are determined by the mass, momentum and energy of the initial datum thanks 
to the conservation properties. As a result of the process of entropy production 
pushing towards local equilibrium combined with the constraints of conservation laws, 
solutions are thus expected to converge to a unique Maxwellian equilibrium.  
Up to a normalization we set without restriction $M(v) = e^{-|v|^2}$ 
as the Maxwellian equilibrium, or equivalently $\rho=\pi^{N/2}$, 
$u=0$ and $T=1/2$. 

The relaxation to equilibrium is studied since the works of Boltzmann and 
it is at the core of the kinetic theory. The motivation is to provide 
an analytic basis for the second principle of thermodynamics for a statistical physics 
model of a gas out of equilibrium. 
Indeed Boltzmann's famous $H$ theorem gives an analytic meaning to 
the entropy production process and identifies possible 
equilibrium states. In this context, proving convergence towards  
equilibrium is a fundamental step to justify Boltzmann model, but 
cannot be fully satisfactory as long as it remains based on non-constructive 
arguments. Indeed, as suggested implicitly by Boltzmann when answering 
critics of his theory based on Poincar\'e recurrence Theorem, 
the validity of the Boltzmann equation breaks for very 
large time (see~\cite[Chapter~1, Section~2.5]{Vi:hb} 
for a discussion). It is therefore crucial to obtain quantitative 
informations on the time scale of the convergence, in order to show that 
this time scale is much smaller than the time scale of 
validity of the model. 
Moreover constructive arguments often provide new qualitative insights into the 
model, for instance here they give a better understanding of 
the dependency of the rate of convergence according to the 
collision kernel and the initial datum. 

\subsection{Assumptions on the collision kernel}

The main physical case of application of this paper is that of 
hard spheres in dimension $N=3$, where (up to a normalization constant)
 \begin{equation}\label{eq:hs}
 B (|v-v_*|, \cos \theta) = |v-v_*|.
 \end{equation}

More generally we shall make the following assumption on the collision kernel:
 \begin{itemize}
 \item[{\bf A.}] We assume that $B$ takes the product form
  \begin{equation}\label{eq:prod}
  B(|v-v_*|,\cos \theta) = \Phi (|v-v_*|) \, b(\cos \theta),
  \end{equation}
 where $\Phi$ and $b$ are nonnegative functions not identically equal to $0$. 
 This decoupling assumption is made for the sake of 
 simplicity and could probably be relaxed at the price of 
 technical complications. 
 \smallskip

 \item[{\bf B.}]  Concerning the kinetic part, we assume 
 $\Phi$ to be given by 
  \begin{equation}\label{eq:hyprad}
  \Phi (z) = C_\Phi \, z^\gamma
  \end{equation}
 with $\gamma \in (0,1]$ and $C_\Phi >0$. 
 It is customary in physics and in mathematics to study the case
 when $\Phi(v-v_*)$ behaves like a power law $|v-v_*|^\gamma$, and one
 traditionally separates between hard potentials ($\gamma >0$),
 Maxwellian potentials ($\gamma = 0$), and soft potentials ($\gamma < 0$). 
 We assume here that we deal with {\bf hard potentials}.
 \smallskip

 \item[{\bf C.}] Concerning the angular part, we assume the control from above
  \begin{equation}\label{eq:grad}
  \forall \, \theta \in [0,\pi], \ \ \ b(\cos \theta) \le C_b.
  \end{equation}
 This assumption is a strong version of Grad's angular cutoff (see~\cite{Grad58}). It is 
 satisfied for the hard spheres model. 

 Moreover, in order to prove the appearance and propagation of exponential moments 
 (see Lemma~\ref{lem:mts}), we shall assume additionally that 
  \begin{equation}\label{eq:hypmts}
  b \mbox{ is nondecreasing and convex on } (-1,1).  
  \end{equation}
 This technical assumption is satisfied for the hard spheres model, since in this 
 case $b$ is constant. 

 Finally, in order to be able to use tools from entropy methods we shall 
 also need the lower bound 
  \begin{equation}\label{eq:beep}
  \forall \, \theta \in [0,\pi], \ \ \ b(\cos \theta) \ge c_b >0.
  \end{equation}
 Again this assumption is trivially satisfied for the hard spheres model. 
 \end{itemize}
 \smallskip

Under these assumptions on the collision kernel $B$, equation~\eqref{eq:base} 
is well-posed in the space of nonnegative solutions with finite and non-increasing 
mass and energy~\cite{MW}. In the sequel by ``solution'' of~\eqref{eq:base} 
we shall always denote these solutions. 

Let us mention that under assumptions \eqref{eq:prod}-\eqref{eq:hyprad}-\eqref{eq:grad}, 
for soft potentials ($\gamma <0$), the linearized operator 
has no spectral gap and no exponential convergence is expected for~\eqref{eq:base} 
(see \cite{Cafl80}). For Maxwellian potentials ($\gamma=0$), exponential convergence 
is known to hold for~\eqref{eq:base} if and only if the initial datum 
has bounded moments of order $s>2$ (see~\cite{CaLu}), and, 
under additionnal moments and smoothness assumptions on the initial datum, 
the rate is known to be governed by the spectral gap of the linearized operator (see \cite{CaGaTo}). 

\subsection{Linearization}

Under assumption~\eqref{eq:grad}, we can define 
  \[ \ell_b := \|b\|_{L^1(\ens{S}^{N-1})} := 
    \big|\ens{S}^{N-2}\big| \, \int_0 ^\pi b(\cos \theta)  \sin^{N-2} \theta \, d\theta < +\infty. \]
Without loss of generality we set $\ell_b =1$ in the sequel. 
Then one can split the collision operator in the following way 
  \begin{eqnarray*}
  Q(g,f) &=& Q^+ (g,f) - Q^- (g,f),  \\
  Q^+ (g,f) &=& \int _{\R^N \times \ens{S}^{N-1}} \Phi(|v-v_*|) \, b(\cos \theta) \, g'_* f' \, dv_* \, d\sigma,  \\
  Q^-(g,f) &=& \int _{\R^N \times \ens{S}^{N-1}} \Phi(|v-v_*|) \, b(\cos \theta) \, g_* f \, dv_* \, d\sigma 
  = (\Phi * g) \, f,
  \end{eqnarray*}
and introduce the so-called {\em collision frequency} 
  \begin{equation} \label{eq:colfreq}
  \nu(v) = \int_{\R^N \times \ens{S}^{N-1}} \Phi(|v-v_*|) \, 
  b(\cos \theta) \, M(v_*) \, dv_* \, d\sigma = (\Phi * M) (v).
  \end{equation}
We denote by $\nu_0 >0$ the minimum value of $\nu$. 

  \begin{definition}[Linearized collision operator]\label{linBo}
  Let $m=m(v)$ be a positive rapidly decaying function. We define 
  the linearized collision operator $\ron{L}_m$ associated with 
  the rescaling $m$, by the formula
    \[ \ron{L}_m (g) = m^{-1} \, \big[ Q(mg,M) + Q(M,mg) \big]. \]
  The particular case when $m=M$ is just called the ``linearized 
  collision operator'', defined by
    \begin{multline*}
    L(h) = M^{-1} \big[ Q(Mh,M) + Q(M,Mh) \big] \\
    = \int_{\R^N \times \ens{S}^{N-1}} \Phi(|v-v_*|) \, b(\cos \theta) \, 
    M(v_*) \big[ h' _*+h' - h_* - h \big] \, dv_* \, d\sigma.
    \end{multline*} 
  \end{definition}

\Remark The linearized collision operator $\ron{L}_m(g)$ corresponds 
to the linearization around $M$ with the scaling $f=M + mg$. Among all 
possible choices of $m$, the case $m=M$ is particular since 
$L$ enjoys a self-adjoint property on the space $g \in L^2(M(v)dv)$, 
which is why this is usually the only case considered. Note that this 
space corresponds to $f \in L^2(M^{-1}(v)dv)$ for the original solution. 
In this paper we shall need other scalings of linearization in order to connect 
the linearized theory to the nonlinear theory. 
We shall use for the scaling function $m(v)$ a ``stretched Maxwellian'' 
of the form $m(v) = \exp \left[ -a |v|^s\right]$ 
with $a >0$ and $0<s<2$ to be chosen later.
\medskip

The linear operators $\ron{L}_m$ splits naturally between 
a multiplicative part $\ron{L}^\nu _m$ and a non-local part $\ron{L}^c _m$ 
(the ``c'' exponent stands for ``compact'' as we shall see) in the following way: 
 \begin{equation}\label{eq:dec1}
 \ron{L}_m(g) = \ron{L}^c _m(g) - \ron{L}^\nu _m(g) \ \ \ \mbox{ with } \ \ \ 
 \ron{L}^\nu _m(g) := \nu \, g
 \end{equation}
where $\nu$ is the collision frequency defined in~\eqref{eq:colfreq},  
and $\ron{L}^c _m$ splits between a ``gain'' part $\ron{L}^+ _m$ (denoted so because it 
corresponds to the linearization of $Q^+$) and a convolution part $\ron{L}^* _m$ as 
 \begin{equation}\label{eq:dec2}
 \ron{L}^c (g) = \ron{L}^+ _m (g) - \ron{L}^* _m (g) 
 \ \ \ \mbox{ with } \ \ \ 
 \ron{L}^* _m (g) :=  m^{-1} \, M \, \big[ (m g) * \Phi \big]
 \end{equation}
and
 \begin{equation}\label{eq:auxdefronL^+}
 \ron{L}^+ _m (g) := m^{-1} \, \int_{\R^N \times \ens{S}^{N-1}} \Phi(|v-v_*|) \, b(\cos \theta) \, 
 \big[ (m g)' M' _*  + M' (m g)' _* \big] \, dv_* \, d\sigma.
 \end{equation}
For $L=\ron{L}_M$ we obtain as a particular case the decomposition
 \begin{equation}\label{eq:dec1:L}
 L(h) = L^c (h) - L^\nu (h) \ \ \ \mbox{ with } \ \ \ 
 L^\nu (h) := \nu \, h
 \end{equation}
and 
 \begin{equation}\label{eq:dec2:L}
 L^c (h) = L^+ (h) - L^* (h)
 \ \ \ \mbox{ with } \ \ \ 
 L^* (h) := ( hM)  * \Phi  
 \end{equation}
and
 \begin{equation}\label{eq:auxdefL^+}
 L^+ (h) := \int_{\R^N \times \ens{S}^{N-1}} \Phi(|v-v_*|) \, b(\cos \theta) \, 
 \big[ h' + h' _* \big] \, M_*  \, dv_* \, d\sigma.
 \end{equation}

\subsection{Spectral theory}

Let us consider a linear unbounded operator $T:\cal{B} \to \cal{B}$ on the 
Banach space $\cal{B}$, defined on a dense domain $\mbox{Dom}(T) \subset \cal{B}$. 
Then we adopt the following notations and definitions:
  \begin{itemize}
  \item we denote by $N(T) \subset \cal{B}$ the {\em null space} of $T$; 
  \item $T$ is said to be {\em closed} if its graph is closed 
  in ${\cal B} \times {\cal B}$;
  \end{itemize}
In the following definitions, $T$ is assumed to be closed.
  \begin{itemize}
  \item the {\em resolvent set} of $T$ denotes 
        the set of complex numbers $\xi$ such that $T-\xi$ is 
        bijective from $\mbox{Dom}(T)$ to $\cal{B}$ and the 
        inverse linear operator $(T-\xi)^{-1}$, 
        defined on $\cal{B}$, is bounded (see \cite[Chapter~3, Section 5]{Kato});
  \item we denote by $\Sigma(T) \subset \C$ the {\em spectrum} of $T$, that is 
        the complementary set of the resolvent set of $T$ in $\C$;
  \item an {\em eigenvalue} is a complex number $\xi \in \C$ such that 
        $N(T-\xi)$ is not reduced to $\{0\}$;
  \item we denote $\Sigma_d(T) \subset \Sigma(T)$ the {\em discrete spectrum} of $T$, 
        {\it i.e.} the set of {\em discrete eigenvalues}, that is the eigenvalues 
        isolated in the spectrum and with finite multiplicity ({\it i.e.} such that the spectral projection  
        associated with this eigenvalue has finite dimension, see \cite[Chapter~3, Section 6]{Kato});
  \item for a given discrete eigenvalue $\xi$, we shall call the {\em eigenspace} 
        of $\xi$ the range of the spectral projection associated with $\xi$;  
  \item we denote $\Sigma_e(T) \subset \Sigma(T)$ the {\em essential spectrum} of $T$ 
        defined by $\Sigma_e(T) = \Sigma(T) \setminus \Sigma_d(T)$;
  \item when $\Sigma(T) \subset \R_-$, we say that $T$ has a {\em spectral gap} when the distance 
  between $0$ and $\Sigma(T) \setminus \{0\}$ is positive, and the spectral gap denotes this distance.
  \end{itemize}
\smallskip

It is well-known from the classical theory of the linearized operator  
(see \cite{Grad63} or \cite[Chapter 7, Section 1]{CIP}) that
 \begin{multline*}
 \langle h,L h \rangle_{L^2(M)} = \int_{\R^N} \bar h \, Lh \, M \, dv = \\ 
 - \frac{1}{4}\int_{\R^N \times 
 \R^N \times \ens{S}^{N-1}} \Phi(|v-v_*|)\, b(\cos \theta) \, 
 \left| h^{'}_{*}+h^{'}-h_{*}-h \right|^2 
 M\, M_{*} \, dv \, dv_{*} \, d\sigma \le 0.
 \end{multline*}
This implies that the spectrum of $L$ in $L^2(M(v)dv)$ is included in $\R_-$. 
Its null space is 
  \begin{equation}\label{noyauL}
  N(L) = \mbox{Span} \left\{ 1, v_1,\dots,v_N, |v|^2 \right\}.
  \end{equation}
These two properties correspond to the linearization of 
Boltzmann's $H$ theorem. 

Let us denote by $D (h)=-\langle h,L h \rangle _{L^2(M)}$ the Dirichlet form for 
$-L$. Since the operator is self-adjoint, the existence of a spectral gap $\lambda$ 
is equivalent to 
 \begin{equation*} \label{eq:sg}
 \forall \, h \bot N(L), \hspace{0.3cm}
          - D(h) \ge \lambda \, \|h\|^2 _{L^2(M)}. 
 \end{equation*}
Controls from below on the collision kernel are necessary so that 
there exists a spectral gap for the linearized operator.  
Concerning the bound from below on $\Phi$, the non-constructive 
proof of Grad suggests that, when the collision kernel 
satisfies Gad's angular cutoff, $L$ has 
a spectral gap if and only if the collision frequency 
is bounded from below by a positive constant ($\nu_0 > 0$). 
Moreover, {\em explicit} estimates on the spectral gap are given in~\cite{BaMo} 
under the assumption that $\Phi$ is bounded from below at infinity, {\it i.e.} 
  \begin{equation*}\label{eq:hardgen}
  \exists \, R \ge 0, \ c_\Phi >0 \ \  ; \ \ 
  \forall \, r \ge R, \ \Phi(r) \ge c_\Phi.
  \end{equation*}
This assumption holds for Maxwellian molecules and hard potentials, with or without angular cutoff. 

Thus under our assumptions on $B$, $L$ has a spectral gap $\lambda \in (0,\nu_0]$ 
(indeed the proof of Grad shows that $\Sigma_e(L)=(-\infty,-\nu_0]$ and the remaining 
part of the spectrum is composed of discrete eigenvalues in $(-\nu_0,0]$ 
since $L$ is self-adjoint nonpositive). Moreover 
as discussed in~\cite[Chapter~4, Section~6]{Ce88}, it was proved in~\cite{KuWi67} 
that $L$ has an infinite number of discrete negative eigenvalues in the interval 
$(-\nu_0,0)$, which implies that 
  \[ 0 < \lambda < \nu_0. \] 
In fact the proof in~\cite{KuWi67} was done for hard spheres, but the 
argument applies to any cutoff hard potential collision kernel as well. 

\subsection{Existing results and difficulties}
 
On the basis of the $H$ theorem and suitable {\em a priori}  
estimates, various authors gave results of 
$L^1$ convergence to equilibrium by compactness arguments 
for the spatially homogeneous Boltzmann equation 
with hard potentials and angular cutoff  
(for instance Carleman~\cite{Carl32}, Arkeryd~\cite{Ar72}, etc.). 
These results provide no information at all on the rate of convergence. 

In~\cite{Grad63} Grad gave the first proof 
of the existence of a spectral gap for the 
linearized collision operator $L$ with hard potentials and 
angular cutoff. His proof was based on 
Weyl's Theorem about compact perturbation and thus did not 
provide an explicit estimate on the spectral gap. 
Following Grad, a lot a works have been done by various authors 
to extend this spectral study to soft potentials (see~\cite{Cafl80}, \cite{GoPo86}), 
or to apply it to the perturbative solutions (see~\cite{Ukai74}) or to 
the hydrodynamical limit (see~\cite{ElPi75} for instance). 

On the basis of these compactness results and linearization tools, 
Arkeryd gave in~\cite{Ar88} the first (non-constructive) proof of 
exponential convergence in $L^1$ for the spatially homogeneous 
Boltzmann equation with hard potentials and angular cutoff. 
His result was generalized to $L^p$ spaces ($1\le p < +\infty$) by Wennberg~\cite{We93}. 
\smallskip

At this point, several difficulties have still to be 
overcome in order to get a quantitative result of exponential convergence:
 \begin{itemize}
 \item[(i)] The spectral gap in $f \in L^2(M^{-1}(v)dv)$ was obtained by non-constructive 
 methods for hard potentials. 
 \item[(ii)] The spectral study was done in the space $f \in L^2(M^{-1}(v)dv)$ 
 for which there is no known {\em a priori} estimate for the nonlinear 
 problem. Matching results obtained in this space and 
 the physical space $L^1((1+|v|^2)dv)$ is one the main difficulties, and 
 was treated in~\cite{Ar88} by a non-constructive argument. 
 \item[(iii)] Finally, any estimate deduced from a linearization 
 argument is valid only in a neighborhood of the equilibrium, and 
 the use of compactness arguments to deduce that the solution 
 enters this neighborhood (as e.g. in~\cite{Ar88}) would 
 prevent any hope of obtaining explicit estimate. 
 \end{itemize} 

First it should be said that in the Maxwellian case, 
all these difficulties have been solved. When the collision kernel 
is independent of the relative velocity,  
Wang-Chang and Uhlenbeck~\cite{WCUh70} and then Bobylev~\cite{Boby88} 
were able to obtain a complete and explicit diagonalization 
of the linearized collision operator, with or without cutoff. Then specific 
metric well suited to the collision operator for Maxwell molecules allowed to 
achieve the goals sketched in the first paragraph of this introduction 
(under additionnal assumptions on the initial datum), see~\cite{CaGaTo} 
and~\cite{CaLu}. However it seems that the proofs of these results 
are strongly restricted to the Maxwellian case. 

In order to solve point~(iii), quantitative estimates in the large 
have been obtained recently, directly on the nonlinear equation,  
by relating the entropy production functional to the relative entropy:
\cite{CC92,CC94,ToVi99,Vi03,MoVi}. 
The latter paper states, for hard potentials with angular cutoff, 
quantitative convergence towards equilibrium 
with rate $O(t^{-\infty})$ for solutions in $L^1((1+|v|^2)dv) \cap L^2$ 
(or only $L^1 ((1+|v|^2)dv)$ in the case of hard spheres). 
However it was proved in~\cite{BobyCerc:conjCer:99} that 
one cannot establish in this functional space a {\em linear} inequality relating 
the entropy production functional and the relative entropy, 
which would yields exponential convergence directly on the nonlinear equation.

Point~(i) was solved in~\cite{BaMo}, which  
gave explicit estimates on the spectral gap for hard potentials, 
with or without cutoff, by relating it explicitly to the one for 
Maxwell molecules. 

In order to solve the remaining obstacle of point~(ii), 
the strategy of this paper is to prove explicit linearized estimates  
of convergence to equilibrium in the space $L^1(\exp \left[ a|v|^s \right]dv)$ with 
$a>0$ and $0<s<\gamma/2$, on which we give explicit results of appearance 
and propagation of the norm, and thus which can be connected to the 
quantitative nonlinear results in~\cite{MoVi}. It will lead us to study 
the linearized operator $\ron{L}_m$ for $m=\exp \left[ -a|v|^s \right]$ on $L^1$, 
which has no hilbertian self-adjointness structure. 

\subsection{Notation} 
In the sequel we shall denote 
$\langle \cdot \rangle = \sqrt{1 + |\cdot|^2}$. 
For any Borel function $w:\R^N \to \R_+$, we define the 
weighted Lebesgue space $L^p (w)$ on $\R^N$ ($p \in [1,+\infty]$), by the norm 
 \[ \| f \|_{L^p(w)} = \left[ \int_{\R^N} |f (v)|^p \, w(v) \, dv \right]^{1/p} \] 
if $p < +\infty$ and 
 \[ \| f \|_{L^\infty (w)} = \sup_{v \in \R^N} |f (v)| \, w(v) \]
when $p = +\infty$. 
The weighted Sobolev spaces $W^{k,p} (w)$ 
($p \in [1,+\infty]$ and $k \in \N$) are defined by the norm 
 \[ \| f \|_{W^{k,p} (w)} = 
       \left[ \sum_{|s| \le k} \|\partial^s f \|^p _{L^p(w)} \right]^{1/p} \]
with the notation $H^k (w) = W^{k,2} (w)$. 
In the sequel we shall denote by $\| \cdot \|$ indifferently 
the norm of an element of a Banach space or the usual 
operator norm on this Banach space, and we shall denote 
by $C$ various positive constants independent of the collision kernel. 

\subsection{Statement of the results}

Our main result of exponential convergence to equilibrium is
  \begin{theorem} \label{theo:cvg}
  Let $B$ be a collision kernel satisfying  
  assumptions~\eqref{eq:prod}, \eqref{eq:hyprad}, \eqref{eq:grad}, \eqref{eq:hypmts}, 
  \eqref{eq:beep}.
  Let $\lambda \in (0,\nu_0)$ be the spectral gap of the linearized operator $L$.   
  Let $f_0$ be a nonnegative initial datum in $L^1(\langle v \rangle^2) \cap L^2$. 
  Then the solution $f(t,v)$ to the spatially 
  homogeneous Boltzmann equation~\eqref{eq:base} with initial datum $f_0$ 
  satisfies: for any $0<\mu \le \lambda$,  
  there is a constant $C >0$, which depends explicitly on $B$, the mass, energy  
  and $L^2$ norm of $f_0$, on $\mu$ and on a lower bound on $\nu_0 -\mu$, such that 
    \begin{equation}\label{eq:cvg}
    \|f(t,\cdot)-M\|_{L^1} \le C \, e^{- \mu t}.
    \end{equation}
  In the important case of hard spheres~\eqref{eq:hs}, the assumption ``$f_0 \in L^1(\langle v \rangle^2) \cap L^2$'' 
  can be relaxed into just ``$f_0 \in L^1 (\langle v \rangle^2)$'', and the same result holds 
  with the constant $C$ in~\eqref{eq:cvg} depending explicitly on $B$, the mass and 
  energy of $f_0$, on $\mu$ and on a lower bound on $\nu_0 -\mu$. 
  \end{theorem}

\Remarks 

1. Note that the optimal rate $\mu=\lambda$ is allowed in the theorem, 
which can be related to the fact that the eigenspace of 
$\ron{L}_m$ associated with the first non-zero eigenvalue $-\lambda$ is not degenerate. 
It seems to be the first time this optimal rate is reached, 
since both the quantitative study in \cite{CaGaTo} 
for Maxwell molecules and the non-constructive results of \cite{Ar88} for hard spheres 
only prove a convergence like $O(e^{-\mu t})$ for any $\mu <\lambda$, where 
$\lambda$ is the corresponding spectral gap. 
\smallskip

2. From~\cite{BaMo}, one deduces the following estimate on $\lambda$: 
when $b$ satisfies the control from below
 \begin{equation*}
 \frac{1}{|\ens{S}^{N-1}|} \, \inf_{\sigma_1, \sigma_2 \in \ens{S}^{N-1}} 
 \int_{\sigma_3 \in \ens{S}^{N-1}} \min \{  b(\sigma_1 \cdot \sigma_3), 
 b(\sigma_2 \cdot \sigma_3) \} \, d\sigma_3 \ge c_b >0
 \end{equation*} 
(which is true for all physical cases, and implied by~\eqref{eq:beep}), then
 \begin{equation*}
 \lambda \ge c_b \, C_\Phi \, 
 \frac{(\gamma/8)^{\gamma/2} \, e^{-\gamma/2} \, \pi}{24}.
 \end{equation*} 
In particular, for hard spheres collision kernels one can compute
 \[ \lambda \ge \pi / (48 \sqrt{2 e}) \approx 0.03. \]  
\medskip

We also state the functional analysis result on the spectrum of $\ron{L}_m$  
used in the proof of Theorem~\ref{theo:cvg} and which has interest in itself. 
We consider the unbounded operator $\ron{L}_m$ on $L^1$ with 
domain $\mbox{Dom}(\ron{L}_m) = L^1 (\langle v \rangle^\gamma)$ and 
the unbounded operator $L$ on $L^2(M)$ with domain 
$\mbox{Dom}(L)= L^2(\langle v \rangle^{2\gamma} M)$. 
These operators are shown to be closed in Proposition~\ref{prop:ferm} 
and Proposition~\ref{prop:ferm:L} and we have
 \begin{theorem}\label{theo:spectre}
 Let $B$ be a collision kernel satisfying  
 assumptions~\eqref{eq:prod}, \eqref{eq:hyprad} and \eqref{eq:grad}. 
 Then the spectrum $\Sigma(\ron{L}_m)$ of $\ron{L}_m$ 
 is equal to the spectrum $\Sigma(L)$ of $L$. 
 Moreover the eigenvectors of $\ron{L}_m$ associated with 
 any discrete eigenvalue are given by those 
 of $L$ associated with the same eigenvalue, multiplied by $m^{-1} M$. 
 \end{theorem}

\Remarks 

1. This theorem essentially means that enlarging the functional space from 
$f \in L^2(M^{-1})$ to $f \in L^1 (m^{-1})$ (for the original solution) 
does not yield new eigenvectors (or additional essential spectrum) 
for the linearized collision operator.
\smallskip

2. It implies in particular that $\ron{L}_m$ only has  
non-degenerate eigenspaces associated with its discrete 
eigenvalues, since this is true for the self-adjoint operator $L$. 
This is related to the fact that the optimal  
convergence rate is exactly $C \, e^{-\lambda t}$ and 
not $C \, t^k \, e^{-\lambda t}$ for some $k >0$. 
It also yields a simple form of the 
first term in the asymptotic developement (see Section~\ref{cvg}).
\smallskip

3. Our study shows that, for hard potentials with cutoff, 
the linear part of the collision operator $f \mapsto Q(M,f) + Q(f,M)$ 
``has a spectral gap'' in $L^1(m^{-1})$, in the sense 
that it satisfies exponential decay estimates 
on its evolution semi-group in this space, with the rate given by the spectral gap of $L$. 
We use this linear feature of the collision process to compensate for 
the fact that the functional inequality ({\em Cercignani's conjecture})
  \begin{equation}\label{eq:CerCon}
  {\cal D}(f) \ge K \, \big[ H(f) - H(M) \big], \quad K>0, 
  \end{equation}
is not true for $f \in L^1(\langle v \rangle^2)$. It also supports 
the fact that \eqref{eq:CerCon} could be true for 
solutions $f$ of~\eqref{eq:base} satisfying some exponential decay at infinity 
(as was questioned in~\cite[Chapter~3, Section~4.2]{Vi:hb}), in the sense $f \in L^1(m^{-1})$. 

\subsection{Method of proof}

The idea of the proof is to establish 
quantitative estimates of exponential decay on the evolution 
semi-group of $\ron{L}_m$. They are used to estimate the rate of convergence 
when the solution is close to equilibrium (where the linear part of the collision 
operator is dominant), whereas the existing nonlinear entropy method, combined  
with some {\em a priori} estimates in $L^1(m^{-1})$, are used to estimate 
the rate of convergence for solutions far from equilibrium. 
The proof splits into several steps.
\smallskip

\noindent 
I. The first step is to prove that $\ron{L}_m$ and $L$ have the same spectrum. 
We use the following strategy: first we localize the essential spectrum of  
$\ron{L}_m$ with the perturbation arguments Grad used for $L$, with 
additional technical difficulties due to the fact that 
the operator $\ron{L}_m$ has no hilbertian self-adjointness structure 
(Proposition~\ref{prop:ess}).   
It is shown to have the same essential spectrum as $L$, 
which is the range of the collision frequency. The main tool is the proof of the 
fact that the non-local part of $\ron{L}_m$ is relatively compact with respect to  
its local part (Lemma~\ref{lem:cpct}). Then, in order to localize the 
discrete spectrum, we show some decay estimates on the 
eigenvectors of $\ron{L}_m$ associated to discrete eigenvalues. 
The operators $\ron{L}_m$ and $L$ are related by 
  \[ \ron{L}_m(g) = m^{-1} M \, L\big( m M^{-1} g \big), \]
and our decay estimates show that any such eigenvector $g$ of $\ron{L}_m$ 
satisfies $m M^{-1} g \in L^2(\langle v \rangle^{2\gamma} M) = \mbox{Dom}(L)$. 
We deduce that $\ron{L}_m$ and $L$ have the same discrete spectrum 
(Proposition~\ref{prop:decvp}). The key tool of the proof is the fact 
that the gain part of $\ron{L}_m$ (as well as the one of $L$) can be  
approximated by some truncation whose range is composed of functions  
with compact support (Propositions~\ref{prop:cvgronL} and~\ref{prop:cvgL}). 
This is where we need the weight $m$ to have exponential decay 
(some polynomial decay would not be sufficient here). 
\smallskip

\noindent
II. The second step is to prove explicit exponential decay estimates 
on the evolution semi-group of $\ron{L}_m$ with optimal rate,  
{\it i.e.} the first non-zero eigenvalue of $\ron{L}_m$ and $L$. To that 
purpose we show sectoriality estimates on $\ron{L}_m$ (Theorem~\ref{theo:sg} and 
Lemma~\ref{lem:resolvent}). 
This requires estimates on the norm of the resolvent of $\ron{L}_m$, 
which are obtained by showing that this norm can be related 
to the norm of the resolvent of $L$ (Proposition~\ref{prop:normR}). 
Again the key tool is the approximation of the gain parts of 
$\ron{L}_m$ and $L$ by some truncation whose range is composed of functions 
with compact support.   
\smallskip

\noindent
III. The third step is the application of these linear estimates 
to the nonlinear problem. A Gronwall argument is used to obtain the 
exponential convergence in an $L^1 (m^{-1})$-neighborhood of the equilibrium 
for the nonlinear problem (Lemma~\ref{lem:Gron}). Moments estimates 
are used to give a new result of appearance and propagation of this 
exponentially weighted norm (Lemma~\ref{lem:mts}), and the nonlinear entropy method 
(in the form of~\cite[Theorems~6.2 and~7.2]{MoVi}) is 
used to estimate the time required to enter this neighborhood (Lemma~\ref{lem:mv}). 

\subsection{Plan of the paper}

Sections~\ref{lin} and~\ref{spec} remain at the functional analysis level. 
In Section~\ref{lin} we introduce suitable approximations of the 
non-local parts of $\ron{L}_m$ and $L$, 
and state and prove various technical estimates on these 
linearized operators useful for the sequel. 
In Section~\ref{spec} we determine the spectrum of $\ron{L}_m$ and show that 
it is equal to the one of $L$. Then in Section~\ref{cvg} we handle solutions of the 
Boltzmann equation: we prove Theorem~\ref{theo:cvg} by translating 
the previous spectral study into explicit estimates on the evolution 
semi-group, and then connecting the latter to the nonlinear theory. 

\section{Properties of the linearized collision operator}
\label{lin}
\setcounter{equation}{0}

In the sequel we fix $m(v) = \exp \left[ -a\, |v|^s \right]$ with $a >0$ and $0<s<2$. The exact values of 
$a$ and $s$ will be chosen later. With no risk of confusion we shall no more write the 
subscript ``$m$'' on the operator $\ron{L}$. We assume in this section 
that the collision kernel $B$ satisfies \eqref{eq:prod}, \eqref{eq:hyprad}, \eqref{eq:grad}. 

\subsection{Introduction of an approximate operator}

Let ${\bf 1}_E$ denote the usual indicator function of the set $E$. 
Roughly speaking we shall truncate smoothly $v$, remove grazing and frontal collisions and 
mollify the angular part of the collision kernel. More precisely, 
let $\Theta: \R \to \R_+$ be an even $C^\infty$ function with 
mass $1$ and support included in $[-1,1]$ and 
$\tilde{\Theta}: \R^N \to \R_+$ a radial $C^\infty$ function 
with mass $1$ and support included in $B(0,1)$. We define the 
following mollification functions ($\epsilon >0$):
  \[ \left\{
     \begin{array}{ll}
     \Theta_\epsilon (x) = \epsilon^{-1} \, \Theta(\epsilon^{-1} x), \ \ \ (x \in \R) \vspace{0.2cm} \\
     \tilde{\Theta}_\epsilon (x) = \epsilon^{-N} \, \tilde{\Theta}(\epsilon^{-1} x), \ \ \ (x \in \R^N).
     \end{array}
     \right. \]
Then for any $\delta \in (0,1)$ we set 
  \begin{equation}\label{eq:auxdefronL^+delta}
  \ron{L}^+ _\delta (g) = {\cal I}_\delta(v) \, m^{-1} \, 
  \int_{\R^N \times \ens{S}^{N-1}} \Phi(|v-v_*|) \, b_\delta (\cos \theta) \, 
  \big[ (m g)' M' _* + M' (m g)' _* \big] \, dv_* \, d\sigma,  
  \end{equation} 
where  
  \begin{equation}\label{eq:defI}
  {\cal I}_\delta = \tilde{\Theta}_\delta * {\bf 1}_{\{|\cdot| \le \delta^{-1}\}}, 
  \end{equation} 
and 
  \begin{equation}\label{eq:defbdelta}
  b_\delta(z) = \left( \Theta_{\delta^2} * {\bf 1}_{\{-1+2\delta^2 \le z \le 1 - 2\delta^2 \}} 
                \right) \, b (z). 
  \end{equation}
We check that 
  \[ \mbox{supp}(b_\delta) \subset \{-1+ \delta^2 \le \cos \theta \le 1 - \delta^2 \}. \]
The approximation induces $\ron{L}_\delta = \ron{L}^+ _\delta - \ron{L}^* - \ron{L}^\nu$, 
following the decomposition (\ref{eq:dec1},\ref{eq:dec2}). 

We define similarly the approximate operator
  \begin{equation}\label{eq:auxdefL^+delta}
  L^+ _\delta (h) = {\cal I}_\delta (v) \, 
  \int_{\R^N \times \ens{S}^{N-1}} b_\delta (\cos \theta) \, \Phi(|v-v_*|) \, 
  \big[ h' + h' _* \big] \, M_* \, dv_* \, d\sigma, 
  \end{equation} 
which induces $L_\delta = L^+ _\delta - L^* - L^\nu$, 
following the decomposition (\ref{eq:dec1:L},\ref{eq:dec2:L}). 

\subsection{Convergence of the approximation}

First for $\ron{L}$ we have 
 \begin{proposition}\label{prop:cvgronL}
 For any $g \in L^1(\langle v\rangle^\gamma)$, we have 
   \[ \left\| \left( \ron{L}^+ - \ron{L}^+ _\delta \right) (g) \right\|_{L^1} \le 
       C_1 (\delta) \, \|g\|_{L^1 (\langle v \rangle^\gamma)} \]
 where $C_1 (\delta) >0$ is an explicit constant depending on the 
 collision kernel and going to $0$ as $\delta$ goes to $0$. 
 \end{proposition}
Before going into the proof of Proposition~\ref{prop:cvgronL}, let us 
enounce as a lemma a simple estimate we shall use several times in the sequel:
  \begin{lemma}\label{lem:mM}
  For all $v \in \R^N$, 
    \begin{equation}\label{mM}
    \left( m M_* (m')^{-1} \right), \ \left( m M_* (m' _*)^{-1}\right)  \, \le \ \ 
    \exp \left[a|v_*|^s - |v_*|^2 \right].
    \end{equation}
  \end{lemma}
\begin{proof}[Proof of Lemma~\ref{lem:mM}]
Indeed 
 \[ \left(m M_* (m')^{-1} \right) = \exp \left[ a |v'|^s - a |v|^s - |v_*|^2 \right] \]
and we have (using the conservation of energy and the fact that $s/2 < \gamma/4 \le 1$)  
 \[ |v'|^s = \left(|v'|^2 \right)^{s/2} \le \left(|v|^2 + |v_*|^2 \right)^{s/2} 
           \le |v|^s + |v_*|^s. \]
This implies immediately~\eqref{mM} (for the other term $\left( m M_* (m' _*)^{-1} \right)$, 
the proof is similar). 
\end{proof}

\begin{proof}[Proof of Proposition~\ref{prop:cvgronL}]
Let us pick $\var >0$. Using the pre-postcollisional change of variable 
\cite[Chapter 1, Section 4.5]{Vi:hb} and the unitary change of 
variable $(v,v_*,\sigma) \to (v_*,v,-\sigma)$, we can write
 \begin{multline*}
 \left\| \left( \ron{L}^+ - \ron{L}^+ _\delta \right) (g) \right\|_{L^1} 
    \le \\
    \int_{\R^N \times \R^N \times \ens{S}^{N-1}} |b-b_\delta| \, |g| \, \langle v \rangle^\gamma \, 
            M_* \, \langle v_* \rangle^\gamma \, m \left[ (m')^{-1} + (m'_*)^{-1} \right] \, dv \, dv_* \, d\sigma \\
        + \int_{\R^N \times \R^N \times \ens{S}^{N-1}} b \, |g| \, \langle v \rangle^\gamma \, 
            M_* \, \langle v_* \rangle^\gamma \, m \left[ (m')^{-1} {\cal C}_\delta(v') 
              + (m'_*)^{-1} {\cal C}_\delta(v'_*) \right] \, dv \, dv_* \, d\sigma \\
          =: I_1 ^\delta + I_2 ^\delta, 
 \end{multline*}
where we have denoted ${\cal C}_\delta (v)= \mbox{Id} - {\cal I}_\delta (v)$, and 
${\cal I}_\delta$ was introduced in \eqref{eq:defI}. 

The goal is to prove that 
  \begin{equation}\label{eq:auxiliaire}
  I_1 ^\delta + I_2 ^\delta \le \var \, \|g\|_{L^1 (\langle v \rangle^\gamma)} 
  \end{equation} 
for $\delta$ small enough. 

By Lemma~\ref{lem:mM},
  \begin{multline*}
  I_1 ^\delta \le 2 \, \int_{\R^N \times \R^N \times \ens{S}^{N-1}} |b-b_\delta| \, |g| \, \langle v \rangle^\gamma \, 
             \langle v_* \rangle^\gamma \, \exp \left[a|v_*|^s - |v_*|^2 \right] 
                           \, dv \, dv_* \, d\sigma \\
             \le 2 \, \|b - b_\delta\|_{L^1(\ens{S}^{N-1})} \, \|g\|_{L^1 (\langle v \rangle^\gamma)} \, 
               \left( \int_{\R^N} \langle v_* \rangle^\gamma \, 
                         \exp \left[a|v_*|^s - |v_*|^2 \right] \, dv_* \right). 
  \end{multline*}
Since $s < \gamma/2 < 2$, we have  
  \[ \left( \int_{\R^N} \langle v_* \rangle^\gamma \, 
             \exp \left[a|v_*|^s - |v_*|^2 \right] \, dv_* \right) <+\infty \]
and thus
  \[ I_1 ^\delta \le C \, \|b - b_\delta\|_{L^1(\ens{S}^{N-1})} \, \|g\|_{L^1 (\langle v \rangle^\gamma)}. \]
Now as $\|b-b_\delta\|_{L^1 (\ens{S}^{N-1})} \to 0$ as $\delta$ goes to $0$ we deduce that there exists 
$\delta_0$ such that for $\delta < \delta_0$ 
  \[ I_1 ^\delta \le (\var/4) \, \|g\|_{L^1 (\langle v \rangle^\gamma)}. \]

For $I_2 ^\delta$, let us denote 
  \[ \phi_1 ^\delta(v) := \int_{\R^N \times \ens{S}^{N-1}} b \,  
            M_* \, \langle v_* \rangle^\gamma \, m \, (m')^{-1} {\cal C}_\delta(v') \, dv_* \, d\sigma \]
  \[ \phi_2 ^\delta(v) := \int_{\R^N \times \ens{S}^{N-1}} b \,  
            M_* \, \langle v_* \rangle^\gamma \, m \, (m'_*)^{-1} {\cal C}_\delta(v'_*) \, dv_* \, d\sigma. \]
We have 
  \[ I_2 ^\delta = \int_{\R^N} |g| \, \langle v \rangle^\gamma \, 
          \left[ \phi_1 ^\delta(v) + \phi_2 ^\delta(v) \right] \, dv. \]
Let us show that $\phi^\delta _1$ and $\phi^\delta _2$ converge to $0$ in $L^\infty$. 
We write the proof for $\phi^\delta _1$, the argument for $\phi^\delta _2$ is symmetric. 
First let us pick $\eta >0$ and introduce the truncation 
$\bar{b}(\cos \theta) = {\bf 1}_{\{-1+\eta \le \cos \theta \le 1 -\eta \}} \, b(\cos \theta)$. 
Then we have 
  \begin{multline*}
  \int_{\R^N \times \ens{S}^{N-1}} |b-\bar{b}| \,  
             M_* \, \langle v_* \rangle^\gamma \, m \, (m'_*)^{-1} {\cal C}_\delta(v'_*) \, dv_* \, d\sigma \\ 
      \le \|b-\bar{b}\|_{L^1 (\ens{S}^{N-1})} \,  
        \left( \int_{\R^N} \langle v_* \rangle^\gamma \, \exp \left[a|v_*|^s - |v_*|^2 \right] \, dv_* \right) 
       \xrightarrow[]{\eta \to 0} 0
  \end{multline*}
and we can choose $\eta$ small enough such that for any $\delta \in (0,1)$
  \[ \left|  \phi_1 ^\delta(v) - \int_{\R^N \times \ens{S}^{N-1}} \bar{b} \,  
             M_* \, \langle v_* \rangle^\gamma \, m \, (m')^{-1} {\cal C}_\delta(v') \, 
             d\sigma \, dv_* \right| \le \frac{\var}8. \]
Second let us pick $R >0$. As 
  \begin{multline*}
  \int_{\{|v_*| \ge R\} \times \ens{S}^{N-1}} \bar{b} \,  
             M_*  \, \langle v_* \rangle^\gamma \, m \, (m')^{-1} {\cal C}_\delta(v') \, d\sigma \, dv_* \\
     \le \|b\|_{L^1 (\ens{S}^{N-1})} \, 
         \left( \int_{\{|v_*| \ge R\}} \langle v_* \rangle^\gamma \, 
         \exp \left[a|v_*|^s - |v_*|^2 \right] \, dv_* \right) 
       \xrightarrow[]{R \to + \infty} 0 
  \end{multline*}
we can choose $R$ large enough such that for any $\delta \in (0,1)$
  \[ \int_{\{|v_*| \ge R\} \times \ens{S}^{N-1}} \bar{b} \,  
             M_* \, \langle v_* \rangle^\gamma \, m \, (m')^{-1} 
                    {\cal C}_\delta(v') \, dv_* \, d\sigma \le \frac{\var}8. \]
Thus we get for any $\delta \in (0,1)$
  \[  \left|  \phi_1 ^\delta(v) - \int_{\{|v_*| \le R\} \times \ens{S}^{N-1}} \bar{b} \,  
             M_* \, \langle v_* \rangle^\gamma \, \
             m \, (m')^{-1} {\cal C}_\delta(v') \, dv_* \, d\sigma \right| \le \frac{\var}4, \]
and it remains to estimate 
  \[ J^\delta (v) = \int_{\{|v_*| \le R\} \times \ens{S}^{N-1}} \bar{b} \,  
             M_* \, \langle v_* \rangle^\gamma \, m \, (m')^{-1} {\cal C}_\delta(v') \, dv_* \, d\sigma. \]
We use the following bound: on the set of angles 
determined by $\{-1+\eta \le \cos \theta \le 1 -\eta \}$, we have 
  \[ |v_* - v'| \le \cos \theta/2 \, |v-v_*| \le \sqrt{1-\eta/2} \, |v-v_*|. \]
Thus when $|v_*| \le R$ we obtain 
  \begin{equation}\label{eq:bornev'}
  |v'| \le R + |v_* - v'| \le R + \sqrt{1-\eta/2} \, |v-v_*| \le 2R + \sqrt{1-\eta/2} \, |v|. 
  \end{equation}
Moreover if we impose $\delta \le (\sqrt{2}R)^{-1}$ (which amounts to take 
a smaller $\delta_0$), we get that if $|v| < \delta^{-1} /\sqrt{2}$ then 
  \[
  |v'| < \sqrt{R^2 + \delta^{-2} /2} \le \delta^{-1}, 
  \]
and thus $J^\delta(v)=0$. So let us assume that $|v| \ge  \delta^{-1} /\sqrt{2}$. For these values of 
$v$ we have (using \eqref{eq:bornev'})
  \[ J^\delta(v) \le \exp \left[ a \big(2R + \sqrt{1-\eta/2} \, |v|\big)^s 
            - a|v|^s \right] \, \|b\|_{L^1 (\ens{S}^{N-1})} \, 
                   \left( \int_{\R^N} M_*  \, \langle v_* \rangle^\gamma \, dv_* \right). \]
To conclude we observe that 
  \[ \exp \left[ a\big(2R + \sqrt{1-\eta/2} \, |v|\big)^s - a|v|^s \right] \xrightarrow[]{|v| \to +\infty} 0 \]
since $s>0$. So, by taking $\delta_0$ small enough, we obtain
  \[ J^\delta(v) \le \frac{\var}8 \]
since it is true for $|v| \ge  \delta^{-1} /\sqrt{2}$ with $\delta$ small 
enough and it is equal to $0$ elsewhere.

This concludes the proof: we have $\| \phi_1 ^\delta \|_{L^\infty} \le 3 \var/8$ 
for $\delta \le \delta_0$. By exactly the same proof we get $\| \phi_2 ^\delta \|_{L^\infty} \le 3 \var/8$, 
and thus $I^\delta _2 \le (3 \var/4) \, \|g\|_{L^1 (\langle v \rangle^\gamma)}$. 
As we had also $I^\delta _1 \le (\var/4) \, \|g\|_{L^1 (\langle v \rangle^\gamma)}$,  
the proof of \eqref{eq:auxiliaire} is complete.  
\end{proof}

\Remark One can see from this proof that we use the fact that the weight 
function $m=m(v)$ satisfies 
  \[ \frac{m(v)}{m(\eta v)} \xrightarrow[]{|v|\to+\infty} 0 \]
for any given $\eta \in [0,1)$. This explains why we do not use a 
polynomial function, but an exponential one. 
\medskip

For $L^+$ we can use classical estimates from Grad \cite{Grad63} to 
obtain a stronger result: the operator $L^+$ is bounded on the space $L^2(M)$, 
and the convergence holds in the sense of operator norm. 

  \begin{proposition}\label{prop:cvgL}
  For any $h \in L^2(M)$, we have 
    \[ \left\| \left( L^+ - L^+ _\delta \right) (h) \right\|_{L^2(M)} \le C_2(\delta) \, \|h\|_{L^2(M)} \]
  where $C_2 (\delta)>0$ is an explicit constant going to $0$ as $\delta$ goes to $0$. 
  \end{proposition}

\begin{proof}[Proof of Proposition~\ref{prop:cvgL}]
Under assumptions~\eqref{eq:hyprad} and~\eqref{eq:grad}, the 
collision kernel $\tilde{B}$ in $\omega$-representation 
\cite[Chapter~1, Section~4.6]{Vi:hb} satisfies
  \[ \tilde{B}(|v-v_*|, \cos \theta) \le 2^{N-2} \, C_b \, C_\Phi \, |v-v_*|^\gamma \, \sin^{N-2} \theta/2 
                            \le 2^{N-2} \, C_b \, C_\Phi \, ( |v-v_*| \sin \theta/2 )^\gamma \]
since $N-2 \ge 1 \ge \gamma$. Hence 
  \[ \tilde{B}(|v-v_*|, \cos \theta) \le 2^{N-2} \, C_b \, C_\Phi \, |v-v'|^\gamma. \]
Then similar computations as in~\cite[Chapter~7, Section~2]{CIP} show that $L^+$ writes
  \[ L^+ (h) (v) = M^{-1/2} (v) \, \int_{u \in \R^N} k(u,v) \, \left( h(u) M^{1/2}(u) \right) \, du \]
with a kernel $k$ satisfying
  \[ k(u,v) \le C \, |u-v|^{1+\gamma-N} \, \exp \left[-\frac{|u-v|^2}{4} - \frac{(|u|^2-|v|^2)^2}{4|u-v|^2}\right]. \]
First we see that this kernel is controlled from above by 
  \[ k(u,v) \le \bar{k}(u-v) := C \, |u-v|^{1+\gamma-N} \, \exp \left[-\frac{|u-v|^2}{4}\right]. \]
Since $\bar{k}$ is integrable on $\R^N$, we deduce that $L^+$ is bounded by  
Young's inequality: 
  \begin{equation}\label{eq:L^+bdd}
  \|L^+\|_{L^2(M)} \le \|\bar{k}\|_{L^1}. 
  \end{equation}

Moreover the computations by Grad~\cite[Section 4]{Grad63} 
(see also~\cite[Chapter 7, Theorem 7.2.3]{CIP}) show that for any $r \ge 0$,
  \[ \int_{\R^N} k(u,v) \, \langle u \rangle^{-r} \, du \le C \, \langle v \rangle^{-r-1} \]
for an explicit constant $C_r>0$. 
Thus if we denote again ${\cal C}_\delta(v)= \mbox{Id} - {\cal I}_\delta(v)$  
(${\cal I}_\delta$ is defined in \eqref{eq:defI}), 
we have for $\|h\|_{L^2(M)} \le 1$: 
  \begin{multline*}
  \big\| {\cal C}_\delta (v) \, L^+(h) \big\|_{L^2 (M)} ^2 
      \le C \, \int_{\R^N} {\cal C}_\delta (v)^2 
         \, \left[ \int_{\R^N} k(u,v) \, M^{1/2}(u) \, h(u) \, du \right]^2 \, dv \\
      \le C \, \int_{\R^N} {\cal C}_\delta (v)^2 \, \left[ \int_{\R^N} k(u,v) \, du \right] 
            \left[ \int_{\R^N} k(u,v) \, M(u) \, h(u)^2 \, du \right] \, dv \\
      \le C \, \int_{\R^N} {\cal C}_\delta (v)^2 \, \langle v \rangle^{-1} \,  
            \left[ \int_{\R^N} k(u,v) \, M(u) \, h(u)^2 \, du \right] \, dv \\
      \le  C \, \langle \delta^{-1} \rangle^{-1} \, 
       \int_{\R^N} \left[ \int_{\R^N} k(u,v) \, M(u) \, h(u)^2 \, du \right] \, dv \\
      \le C \, \langle \delta^{-1} \rangle^{-1} \, 
      \left[ \int_{\R^N} \left( \int_{\R^N} k(u,v) \, dv \right) \, M(u) \, h(u)^2 \, du \right] \\
      \le C \, \langle \delta^{-1} \rangle^{-1} \, 
           \left[ \int_{\R^N} \, M(u) \, h(u)^2 \, du \right] \le C \, \langle \delta^{-1} \rangle^{-1} 
  \end{multline*}
using finally the $L^1$ bound 
  \[ \int_{\R^N} k(u,v) \, dv \le \int_{\R^N} \bar{k}(u-v) \, dv \le \|\bar{k}\|_{L^1} < +\infty \]
independent of $u$. 
This shows that 
  \begin{equation}\label{eq:auxdecgrad}
  \big\| {\cal C}_\delta (v) \, L^+ \big\|_{L^2 (M)} = O(\delta^{1/2}) 
  \end{equation}
and thus ${\cal C}_\delta (v) L^+$ goes to $0$ as $\delta$ goes to $0$ in the sense of operator 
norm, with explicit rate.  

Let us again pick $h$ with $\|h\|_{L^2(M)} \le 1$, then 
  \begin{equation}\label{eq:auxtriang}
  \left\| \left( L^+ - L^+ _\delta \right) (h) \right\|_{L^2(M)} \le 
     \| {\cal C}_\delta (v) \, L^+(h) \|_{L^2 (M)} + 
        \| {\cal I}_\delta (v) \, L^+ _{|b-b_\delta|} (h) \|_{L^2(M)} 
  \end{equation}
where the notation $L^+ _{|b-b_\delta|}$ stands for the linearized collision operator $L^+$ 
with the collision kernel $\Phi \, |b-b_\delta|$ instead of $\Phi \, b$. 
We have 
  \begin{multline*} 
  \big\| {\cal I}_\delta (v) \, L^+ _{|b-b_\delta|} (h) \big\|^2 _{L^2(M)} \\
      \le C \, \int_{\R^N} {\cal I}_\delta (v) \, \langle v \rangle^{2\gamma}
        \left[ \int_{\R^N \times \ens{S}^{N-1}} M_* \, |b-b_\delta| \, 
         \langle v_* \rangle^{\gamma} \, \left(h' + h' _* \right) \, dv_* \, d\sigma \right]^2 \, M(v) \, dv. 
  \end{multline*} 
We use the truncation of $v$ to control $\langle v \rangle^{\gamma}$ and the Cauchy-Schwarz 
inequality together with the bound 
  \[ \int_{\R^N} M(v_*) \, \langle v_* \rangle^{2\gamma} \, dv_* < +\infty. \]
This yields
  \begin{multline*}
  \big\| {\cal I}_\delta (v) \, L^+ _{|b-b_\delta|} (h) \big\|^2 _{L^2(M)} \\ 
      \le C \, \langle \delta^{-1} \rangle^{2\gamma} \, \|b-b_\delta\|_{L^1(\ens{S}^{N-1})} \, 
        \int_{\R^N \times \R^N \times \ens{S}^{N-1}}  |b-b_\delta| \, 
           \left[ (h')^2 + (h' _*)^2 \right] \, M \, M_* \, dv \, dv_* \, d\sigma. 
  \end{multline*}
Then using the pre-postcollisional change of variable we get 
 \begin{multline*}
 \big\| {\cal I}_\delta(v) \, L^+ _{|b-b_\delta|} (h) \big\|^2 _{L^2(M)} \\
 \le C \, \langle \delta^{-1} \rangle^{2\gamma} \, \|b-b_\delta\|_{L^1(\ens{S}^{N-1})} \, 
        \int_{\R^N \times \R^N \times \ens{S}^{N-1}}  |b-b_\delta| \, 
           \left[ h^2 + (h _*)^2 \right] \, M \, M_* \, dv \, dv_* \, d\sigma \\
     \le C \, \langle \delta^{-1} \rangle^{2\gamma} \, 
      \|b-b_\delta\|_{L^1(\ens{S}^{N-1})}^2. 
 \end{multline*}
Finally by \eqref{eq:defbdelta} we have 
  \[ \|b-b_\delta\|_{L^1(\ens{S}^{N-1})} ^2 \le C \, \delta^4 \]
and we deduce 
  \[  \big\| {\cal I}_\delta(v) \, L^+ _{|b-b_\delta|} (h) \big\|^2 _{L^2(M)} 
     \le C \, \langle \delta^{-1} \rangle^{2\gamma} \, \delta^4. \]
Since $\gamma \le 1$ and we have 
  \[ \big\|{\cal I}_\delta (v) \, L^+ _{|b-b_\delta|}\big\|_{L^2(M)} = O(\delta^{2-\gamma}), \]
we deduce that ${\cal I}_\delta (v) L^+ _{|b-b_\delta|}$ goes to $0$ as $\delta$ goes to $0$ 
in the sense of operator norm, with explicit rate. 
Together with \eqref{eq:auxdecgrad} and \eqref{eq:auxtriang}, this concludes the proof. 
\end{proof}

\subsection{Estimates on $\ron{L}$}

  \begin{proposition}\label{prop:ronL}
  For any $\delta \in (0,1)$, we have the following properties.  
    \begin{enumerate}
    \item[(i)] There exists an explicit constant $C_3 >0$ depending  
    only on the collision kernel such that 
      \begin{equation}\label{eq:gain:poids}
      \left\{
       \begin{array}{ll}
       \| \ron{L} ^+ (g) \|_{L^1} \le& C_3 \, 
               \| g \|_{L^1 (\langle v \rangle^\gamma)} \vspace{0.2cm} \\
       \| \ron{L} ^+ _\delta (g) \|_{L^1} \le& 
              C_3 \, \| g \|_{L^1 (\langle v \rangle^\gamma)}.  
       \end{array}
      \right.  
      \end{equation}
    \item[(ii)] 
    There exists an explicit constant $C_4 (\delta) >0$ 
    depending on $\delta$ (and going to infinity as $\delta$ goes to $0$) such that 
      \begin{equation}\label{eq:dec:gain}
      \forall \, v \in \R^N, \ \ \ 
      \left| \ron{L}^+ _\delta (g) (v) \right|
      \le C_4 (\delta) \, {\cal I}_\delta (v) \, \| g \|_{L^1}.
      \end{equation} 
    \item[(iii)] 
    There is an explicit constant $C_5 (\delta) >0$ depending on $\delta$ (and possibly 
    going to infinity as $\delta$ goes to $0$) 
    such that for all $\delta \in (0,1)$
      \begin{equation}\label{eq:reg:gain}
      \|\ron{L}^+ _\delta (g)\|_{W^{1,1}}  \le C_5 (\delta) \, \| g \|_{L^1}.
      \end{equation}
    \item[(iv)] There is an explicit constant $C_6 > 0$ such that
      \begin{equation}\label{eq:dec:M}
      \forall \, v \in \R^N, \ \ \ 
         |\ron{L}^* (g) (v)| \le 
      C_6 \, \| g \|_{L^1} \, m^{-1} \, \langle v \rangle^\gamma \, M(v).  
      \end{equation}
    \item[(v)] 
    There exists an explicit constant $C_7 >0$ such that 
      \begin{equation}\label{eq:reg:M}
      \|\ron{L}^* (g)\|_{W^{1,1}}  \le C _7 \, \| g \|_{L^1}.  
      \end{equation}
    \item[(vi)] There are some explicit constants $n_0$, $n_1>0$ such that
      \begin{equation}\label{eq:nu}
      \forall \, v \in \R^N, \ \ \ 
      n_0 \, \langle v \rangle^\gamma \le \ 
      \nu(v) \ \le n_1 \, \langle v \rangle^\gamma.
      \end{equation}
    \end{enumerate}
  \end{proposition}
\Remark
The regularity property~\eqref{eq:reg:gain} is proved 
here by direct analytic computations on the kernel but is  
reminiscient of the regularity property on $Q^+$ of the form  
 \[ \| Q^+ (g,f) \|_{H^s} \le C \| g \|_{L^1 _2} \| f \|_{L^2 _\gamma} \]  
with $s>0$ (see~\cite{Lions94,We94,BoDe98,Lu98,MoVi}), 
proved with the help of tools from harmonic analysis to handle 
integral over moving hypersurfaces. Here we do not need such tools since 
the function integrated on the moving hyperplan is just a gaussian. 
Note that, using the arguments from point~(iii), one  
could easily prove that $L^+ _\delta$ is bounded from $L^2(M)$ into 
$H^1(M)$. Since $L^+ _\delta$ converges to $L^+$, 
this would provide a proof for the compactness 
of $L^+$, alternative to the one of Grad (which is based on the Hilbert-Schmidt theory). 
But in fact most of the key estimates of the proof of 
Grad were used in the proof of Proposition~\ref{prop:cvgL}. 
Nevertheless it underlines the fact that the compactness property of $L^+$ 
can be linked to the same kind of regularity effect that we observe for 
the nonlinear operator $Q^+$. 
\medskip

\begin{proof}[Proof of Proposition~\ref{prop:ronL}]
Point~(i) follows directly from convolution-like estimates in~\cite[Section~2]{MoVi} 
together with inequality~\eqref{mM}. Point (ii) is a direct consequence of the estimates 
  \[ Q^+: L^\infty (\langle v \rangle^\gamma) \times L^1 (\langle v \rangle^\gamma) \to L^\infty \]
  \[ Q^+: L^1 (\langle v \rangle^\gamma) \times L^\infty (\langle v \rangle^\gamma) \to L^\infty \]
valid when grazing and frontal collisions are removed (see~\cite[Section~2]{MoVi} again) 
and thus valid for the quantity
  \[  m^{-1} \, \int_{\R^N \times \ens{S}^{N-1}} \Phi(|v-v_*|) \, b_\delta (\cos \theta) \, 
  \big[ (m g)' M' _* + M' (m g)' _* \big] \, dv_* \, d\sigma \]
appearing in the formula of $\ron{L}^+ _\delta$.  
Point~(iv) is trivial and point~(vi) is well-known. It remains to prove the regularity estimates. 

For the regularity of $\ron{L}^+ _\delta$ (point~(iii)), 
we first derive a representation in the spirit of the computations 
of Grad (it is also related to the Carleman representation, see~\cite[Chapter~1, Section~4.6]{Vi:hb}).  
Write the collision integral with the ``$\omega$-representation'' 
(see~\cite[Chapter~1, Section~4.6]{Vi:hb} again)
  \begin{multline}\label{eq:auxomega}
  \ron{L}^+ _\delta (g) = {\cal I}(v) \, 
  m^{-1} \, \int_{\R^N \times \ens{S}^{N-1}} \Phi(|v-v_*|) \, 
  \tilde{b}_\delta \left(\omega \cdot \frac{v-v_*}{|v-v_*|}\right) \\ 
  \big[ (m g)' M' _* + M' (m g)' _* \big] \, d\omega \, dv_{*}.
  \end{multline}
In this new parametrization of the collision, the velocities 
before and after collision are related by 
  \[ v' = v + ((v_*-v) \cdot \omega) \, \omega, \qquad
     v'_* = v_* - ((v_*-v) \cdot \omega) \, \omega \]
and the angular collision kernel is given by 
  \[ \tilde{b}_\delta(u) = 2^{N-1} \, u^{N-2} \, b_\delta(1-2u^2). \]
Modulo replacing $\tilde{b}_\delta$ by a symmetrized version 
$\tilde{b}^s _\delta (\theta) = \tilde{b}_\delta (\theta) 
+ \tilde{b}_\delta (\pi/2 -\theta)$, we can combine terms appearing 
in \eqref{eq:auxomega} into just one: 
  \begin{equation*}
  \ron{L}^+ _\delta (g) = {\cal I}(v) \, m^{-1} \, \int_{\R^N \times \ens{S}^{N-1}} 
  \tilde{b}^s _\delta \left(\omega \cdot \frac{v-v_*}{|v-v_*|}\right) \, \Phi(|v-v_*|) \, 
  (m g)' M' _* \, dv_* \, d\omega.
  \end{equation*}
Then, keeping $\omega$ unchanged, we make the translation change of variable 
$v_* \to V = v_* - v$, 
  \begin{multline*}
  \ron{L}^+ _\delta (g)(v) = {\cal I}(v) \, m^{-1} \, \int_{\R^N \times \ens{S}^{N-1}} 
  \tilde{b}^s _\delta \left(\omega \cdot \frac{V}{|V|} \right) \, \Phi(|V|) \\  
  (m g)\big(v + (V \cdot \omega) \omega\big)  \, M( v + V + (V \cdot \omega) \omega) \, dV \, d\omega.
  \end{multline*}
Then, still keeping $\omega$ unchanged, we write the orthogonal decomposition 
$V = V_1 \omega + V_2$ with $V_1 \in \R$ and 
$V_2 \in \omega^\bot$ (the latter set can be identified with $\R^{N-1}$) 
  \begin{multline*}
  \ron{L}^+ _\delta (g)(v) = {\cal I}(v) \, m^{-1} \, \int_{\ens{S}^{N-1} \times \R \times \omega^\bot}
  \tilde{b}^s _\delta \left(\frac{V_1}{\sqrt{|V_1|^2+|V_2|^2}} \right) \\ 
  \Phi\big(\sqrt{|V_1|^2+|V_2|^2}\big) \, \, 
  (m g)\big(v + V_1 \omega\big) \, M( v + V_2) \, d\omega \, dV_1 \, dV_2.
  \end{multline*}
Finally we reconstruct the polar variable $W = V_1 \, \omega$ 
  \begin{multline*}
  \ron{L}^+ _\delta (g)(v) = {\cal I}(v) \, m^{-1} \, \int_{\R^N \times W^\bot} |W|^{-(N-1)} \, 
  \tilde{b}^s _\delta \left(\frac{|W|}{\sqrt{|W|^2+|V_2|^2}} \right) \\ \Phi\big(\sqrt{|W|^2+|V_2|^2}\big) \, 
  \, (m g)\big(v + W\big)  \, M( v + V_2) \, dW \, dV_2.
  \end{multline*}
This finally leads to the following representation of $\ron{L}^+ _\delta$:
  \begin{multline}\label{eq:alter}
  \ron{L}^+ _\delta (g)(v) = {\cal I}(v) \, m^{-1} \, \int_{\R^N}  (m g)\big(v + W\big) \times \\ 
  \Bigg( \int_{W^\bot} |W|^{-(N-1)} \, \tilde{b}^s _\delta \left(\frac{|W|}{\sqrt{|W|^2+|V_2|^2}} \right) \, 
  \Phi\big(\sqrt{|W|^2+|V_2|^2}\big) \, M( v + V_2) \, dV_2 \Bigg) \, dW.
  \end{multline}

Then we compute a derivative along some coordinate $v_i$. By integration by parts,
  \begin{multline*}
  \partial_{v_i} \ron{L}^+ _{\delta} (g)(v) =   
  - {\cal I}(v) \, m^{-1} \, \int_{\R^N}  (m g)(v + W) \times \\
  \partial_{W_i} \left[  \int_{W^\bot}  |W|^{-(N-1)} \, \tilde{b}^s _\delta \,   
  \left(\frac{|W|}{\sqrt{|W|^2+|V_2|^2}} \right) \, 
  \Phi \big(\sqrt{|W|^2+|V_2|^2}\big) 
  \, M( v + V_2) \, dV_2 \right] \, dW \\
  +  {\cal I}(v) \, m^{-1} \, \int_{\R^N}  (m g)(v + W) \times \\
   \left[ \int_{W^\bot} |W|^{-(N-1)} \, \tilde{b}^s _\delta \, \left(\frac{|W|}{\sqrt{|W|^2+|V_2|^2}} \right) \, 
  \Phi \big(\sqrt{|W|^2+|V_2|^2}\big) 
  \, \partial_{v_i} M( v + V_2) \, dV_2 \right] \, dW \\
  +  \partial_{v_i} ( {\cal I}(v) m^{-1}) \, \int_{\R^N}  (m g)(v + W) \times \\
  \left[ \int_{W^\bot} |W|^{-(N-1)} \, \tilde{b}^s _\delta \, \left(\frac{|W|}{\sqrt{|W|^2+|V_2|^2}} \right) \, 
  \Phi \big(\sqrt{|W|^2+|V_2|^2}\big) 
  \,  M( v + V_2) \, dV_2 \right] \, dW \\
  =: I_1 + I_2 + I_3. 
  \end{multline*}
The functions ${\cal I}_\delta m^{-1}$ and $\partial_{v_i} ({\cal I}_\delta (v) m^{-1})$ 
are bounded in the domain of 
truncation. Concerning the term $I_2$ we have immediately  
  \[ |\partial_{v_i} M| 
     \le C \, M^{1/2} \] 
and thus straightforwardly 
  \[ \int_{\R^N} I_2 \, dv, \ \ \int_{\R^N} I_3 \, dv  \le 
     C(\delta) \, \|mg\|_{L^1 (\langle v \rangle^\gamma)} \le C(\delta) \, \|g\|_{L^1}. \]

For the term $I_1$, we use the fact that, in the domain of the angular truncation $b_\delta$, we have 
  \begin{equation}\label{eq:WV_2}
  \alpha_\delta |V_2| \le |W| \le \beta_\delta |V_2| 
  \end{equation}
for some constants $\alpha_\delta>0$ and 
$\beta_\delta>0$ depending on $\delta$. In order not to 
deal with a moving domain of integration we shall write the integral 
as follows. Since the integral is even with respect to $W$, we can restrict the study 
to the set of $W$ such that the first coordinate $W_1$ is nonnegative. 
We denote $e_1$ the first unit vector of the corresponding orthonormal basis. 
Then we define the following orthogonal linear transformation of $\R^N$, for some 
$\omega \in \ens{S}^{N-1}$: 
  \[ \forall \, X \in \R^N, \ \ \ R(\omega,X) = 2 \frac{(e_1 + \omega) \cdot X}{|e_1+\omega|^2} (e_1 + \omega) - X. \]
Geometrically $R(\omega,\cdot)$ is the axial symmetry with respect to the line defined by the vector $e_1 + \omega$. 
It is straightforward that $R(\omega,\cdot)$ is a unitary diffeomorphism from $\{X, \ X_1=0\}$ onto 
$\omega^\bot$. We deduce that 
  \begin{multline*}
  \int_{W^\bot}  |W|^{-(N-1)} \, \tilde{b}^s _\delta \,   
  \left(\frac{|W|}{\sqrt{|W|^2+|V_2|^2}} \right) \, 
  \Phi \big(\sqrt{|W|^2+|V_2|^2}\big) \, M( v + V_2) \, dV_2 \\
  = \int_{\R^{N-1}} |W|^{-(N-1)} \, \tilde{b}^s _\delta \,  
  \left(\frac{|W|}{\sqrt{|W|^2+|U|^2}} \right) \\  
  \Phi \big(\sqrt{|W|^2+|U|^2}\big) \, M\left[ v + R\left(\frac{W}{|W|},(0,U)\right)\right]  \, dU.  
  \end{multline*}

Thus we compute by differentiating each term
  \begin{multline*}
  \partial_{W_i} \left[ |W|^{-(N-1)} \, \tilde{b}^s _\delta \, 
  \left(\frac{|W|}{\sqrt{|W|^2+|U|^2}}\right) \, 
  \Phi \big(\sqrt{|W|^2+|U|^2}\big) \, M\left[ v + R\left(\frac{W}{|W|},(0,U)\right)\right] \right] \\
  = \left[ -(N-1) \, |W|^{-N} \, \frac{W_i}{|W|} \, \tilde{b}^s _\delta 
  \, \Phi \, M \right] \\
  + \left[  |W|^{-(N-1)} \, \frac{W_i |U|^2}{|W|(|W|^2+|U|^2)^{3/2}} \,  
  (\tilde{b}^s _{\delta})' \,  \Phi \, M \right] \\
  + \left[  |W|^{-(N-1)} \, \frac{W_i}{\sqrt{|W|^2+|U|^2}} \, 
  \tilde{b}^s _\delta \, \Phi' \, M \right] \\
  +  \left[ - |W|^{-(N-1)} \, \partial_{W_i} \left|v+R\left(\frac{W}{|W|},(0,U)\right)\right|^2 \, 
  M \, \tilde{b}^s _\delta \, \Phi \right] \\
  =: I_{1,1} + I_{1,2} + I_{1,3} + I_{1,4}.
  \end{multline*}
Then we use the fact that $\tilde{b}^s _\delta$ and $(\tilde{b}^s _\delta)'$ 
are bounded in $L^\infty$ by some constant depending on $\delta$, and 
  \[ \Phi (z) \le C_\Phi \, z^\gamma, \hspace{0.3cm} |\Phi' (z)| \le C_\Phi \, \gamma \, z^{\gamma-1}. \] 
We write the previous expression according to $|W|$ only thanks 
to~\eqref{eq:WV_2} (using that $|U|=|V_2|$). The three first terms are controlled as follows 
  \[ I_{1,1}, I_{1,2} \le C(\delta) \, |W|^{-N} \, 
  \Phi \big(\sqrt{|W|^2+|U|^2}\big) \, M \le C(\delta) \, |W|^{-N+\gamma} \, M \]
  \[ I_{1,3} \le |W|^{-N+1} \, 
  \big|\Phi' \big(\sqrt{|W|^2+|U|^2}\big)\big| \, M \le C(\delta) \, |W|^{-N+\gamma} \, M, \]
for some constant $C(\delta)$ depending on $\delta$. 
Finally for the fourth term, easy computations give
  \[ \partial_{W_i} \left|v+R\left(\frac{W}{|W|},(0,U)\right)\right|^2 \le C \, \left( \frac{1+|v|^2+|U|^2}{|W|} \right) \]
and thus using the controls~\eqref{eq:WV_2} we deduce that on the domain 
of truncation for $v$ we have  
  \[ I_{1,4} \le C(\delta) \, \left( |W|^{-N+\gamma} + |W|^{-N+1+\gamma} \right) \, M. \]

Thus $I_1$ is controlled by 
  \begin{multline*}
  I_1 \le C(\delta) \, {\cal I}_\delta (v) \, m^{-1}(v) \, 
  \int_{\R^N} (m |g|)(v+W) \, \left( |W|^{-N+\gamma} + |W|^{-N+1+\gamma} \right) \\ \left[ \int_{\R^{N-1}} 
  M \left(v+R\left(\frac{W}{|W|},(0,U)\right)\right) \, dU \right] \,dW \\
  = C(\delta) \, {\cal I}_\delta (v)  \, m^{-1}(v) \, 
  \int_{\R^N} (m |g|)(v+W) \, \left( |W|^{-N+\gamma} + |W|^{-N+1+\gamma} \right) \, M(v\cdot W/|W|) \,dW 
  \end{multline*}
for some new constant $C(\delta)$. Hence, 
using that $|\cdot|^{-N+\gamma}$ and $|\cdot|^{-N+1+\gamma}$ 
are integrable near $0$ in $\R^N$ (as $\gamma >0$) 
and a translation change of variable $u=v+W$, we obtain
  \begin{multline*}
  \int_{\R^N} I_1 \, dv \\ \le C(\delta) \, 
  \int_{\R^N \times \R^N} {\cal I}_\delta (v)  \, m^{-1}(v) \, (m |g|)(v+W) \, 
  \left( |W|^{-N+\gamma} + |W|^{-N+1+\gamma} \right) \, 
  M(v\cdot W/|W|) \, dv \, dW \\
  \le C(\delta) \, \int_{\R^N} (m |g|)(u) \, \left[ \int_{\R^N} {\cal I}_\delta (u-W)  \, m^{-1}(u-W) \, 
  \left( |W|^{-N+\gamma} + |W|^{-N+1+\gamma} \right) \, dW \right] \, du \\
  \le C(\delta) \, \| g\|_{L^1}
  \end{multline*}
for some new constant $C(\delta)$ (for the last inequality 
we use the fact that the truncation ${\cal I}_\delta$ reduces the integration 
over $W$ to a bounded domain). We deduce that 
 \[ \int_{\R^N} I_1 \, dv  \le C(\delta) \, \| g \|_{L^1}. \]

Gathering the estimates for $I_1$, $I_2$, $I_3$, we obtain 
 \[  \| \partial_{v_i} \ron{L}^+ _\delta (g) \|_{L^1} 
 \le C(\delta) \, \| g \|_{L^1}. \]
The proof of inequality~\eqref{eq:reg:gain} 
is completed by writing this estimate on each derivative and 
using the bound~\eqref{eq:dec:gain} on $\|\ron{L}^+ _\delta\|_{L^1}$. 

Finally point (v) is simpler since $\ron{L}^*$ has a more classical convolution structure. 
We compute a derivative along some coordinate $v_i$:
  \begin{multline*}
  \partial_{v_i} \ron{L}^* (g)(v) =  
  \left( \int_{\R^N} 
  \partial_{v_i} \Phi (v-v_*) (m g)(v_*) \, dv_* \right) 
  \, m^{-1}(v) \, M(v) \\
  + \left( \int_{\R^N} \Phi (v-v_*) (m g)(v_*) \, dv_* \right) \, 
     \partial_{v_i} \left(m^{-1}(v) \, M(v) \right)
  \end{multline*}
and it is straightforward to control the $L^1$ norm of these two terms 
according to the $L^1$ norm of $g$. 
\end{proof}

Now we can deduce some convenient properties in order to handle the 
operator $\ron{L}$ with the tools from the spectral theory. We give properties for each part of 
the decomposition as well as for the global operator. 
  \begin{proposition}\label{prop:ferm} 
    \begin{enumerate}
    \item[(i)] For all $\delta \in (0,1)$, 
    the operator $\ron{L}^+ _\delta$ is bounded on $L^1$ (with 
    explicit bound $C_5 (\delta)$). The operator $\ron{L}^+$, 
    with domain $L^1 (\langle v \rangle^\gamma)$, 
    is closable on $L^1$ . 
    \item[(ii)] The operator $\ron{L}^*$ is bounded 
    on $L^1$ (with explicit bound $C_7$).
    \item[(iii)] The operator $\ron{L}^\nu$, with domain  
    $L^1 (\langle v \rangle^\gamma)$, is closed on $L^1$ . 
    \item[(iv)] 
    The operator $\ron{L}$, with domain $L^1 (\langle v \rangle^\gamma)$, 
    is closed on $L^1$.   
    \end{enumerate} 
  \end{proposition}
\begin{proof}[Proof of Proposition~\ref{prop:ferm}]
For point~(i), the boundedness of 
$\ron{L}^+ _\delta$ is already proved in~\eqref{eq:reg:gain} 
and $\ron{L}^+$ is well-defined on $L^1 (\langle v \rangle^\gamma)$ 
from~\eqref{eq:gain:poids}. 
Let us prove that $\ron{L}^+$ is closable when defined on $L^1$ with domain 
$L^1 (\langle v \rangle^\gamma)$. It means that 
for any sequence $(g_n)_{n\ge0}$ in $L^1 (\langle v \rangle^\gamma)$, 
going to $0$ in $L^1$, and such that 
$\ron{L}^+ (g_n)$ converges to $G$ in $L^1$, 
we have $G \equiv 0$. We can write 
  \[ \ron{L}^+ (g_n) = m^{-1} \, \ron{\bar{L}}^+ (g_n) \]
where 
  \[ \ron{\bar{L}}^+ (g_n) = Q^+ (M, m g_n) + Q^+ (m g_n, M). \] 
It is straightforward to see from the proof of~\eqref{eq:gain:poids} 
(using~\cite[Theorem~2.1]{MoVi}) that 
$\ron{\bar{L}}^+$ is bounded in $L^1$. So 
$g_n \to 0$ in $L^1$ implies that 
$\ron{\bar{L}}^+(g_n) \to 0$ in $L^1$, which 
implies that $\ron{\bar{L}}^+(g_n)$ 
goes to $0$ almost everywhere, up to an extraction. 
After multiplication by 
$m^{-1}$, we deduce that, up to an extraction, 
$\ron{L}^+(g_n)$ goes to $0$ almost everywhere. This implies that 
$G \equiv 0$ and concludes the proof. 

Point~(ii) is already proved in~\eqref{eq:reg:M}. For point~(iii), $\ron{L}^\nu$ 
is well-defined on $L^1 (\langle v \rangle^\gamma)$ from~\eqref{eq:nu} and the closure 
property is immediate: for any sequence $(g_n)_{n\ge0}$ 
in $L^1 (\langle v \rangle^\gamma)$ such that $g_n \to g$ in $L^1$ 
and $\ron{L}^\nu (g_n) \to G$ in $L^1$, 
we have, up to an extraction, that $g_n$ goes to $g$ almost everywhere 
and $\nu g_n$ goes to $G$ almost 
everywhere. So $G = \nu g = \ron{L}^\nu (g)$ almost everywhere, 
and moreover as $G \in L^1$, we deduce from~\eqref{eq:nu} 
that $g \in L^1 (\langle v \rangle^\gamma)=\mbox{Dom}(\ron{L}^\nu)$.

For point~(iv), first we remark that $\ron{L} _\delta$ is trivially closed since it is the 
sum of a closed operator plus a bounded operator (see~\cite[Chapter~3, Section~5.2]{Kato}). 
In order to prove that $\ron{L}$ is closed we shall prove 
on $\ron{L}^c$ a quantitative relative 
compactness estimate with respect to $\nu$. 

By Proposition~\ref{prop:cvgronL} we have 
  \[ \|\ron{L}^+ (g) \|_{L^1} \le \|\ron{L}^+ _\delta (g) \|_{L^1} 
                                 + \frac{n_0}{2} \, \|g\|_{L^1 (\langle v \rangle^\gamma)} \] 
if we choose $\delta > 0$ such that $C_1 (\delta) \le n_0 /2$ ($n_0 >0$ is defined 
in \eqref{eq:nu}). Hence
  \begin{equation}\label{eq:auxcomp}
  \|\ron{L}^+ (g) \|_{L^1} \le C_4(\delta) \, \| g \|_{L^1} 
                                 + \frac{n_0}{2} \, \|g\|_{L^1 (\langle v \rangle^\gamma)} 
  \end{equation}
and thus for the whole non-local part 
  \[ \|\ron{L}^c (g) \|_{L^1} \le \left[ C_4 (\delta) + C_7 \right] \, \| g \|_{L^1} 
               + \frac{n_0}{2} \, \|g\|_{L^1 (\langle v \rangle^\gamma)}
            = C \, \| g \|_{L^1} + \frac{n_0}{2} \, \|g\|_{L^1 (\langle v \rangle^\gamma)}. \]
Then by triangular inequality 
  \[  \|\ron{L} (g) \|_{L^1} 
  \ge \|\ron{L}^\nu (g) \|_{L^1} - \|\ron{L}^c (g) \|_{L^1} 
  \ge \frac{n_0}{2} \, \|g\|_{L^1 (\langle v \rangle^\gamma)} - C \, \| g \|_{L^1} \]
which implies that 
  \begin{equation}\label{eq:compprec}
  \|g\|_{L^1 (\langle v \rangle^\gamma)} 
  \le \frac{2}{n_0} \, \Big( \|\ron{L} (g) \|_{L^1} 
    + C \, \| g \|_{L^1} \Big).
  \end{equation}
If $(g_n)_{n \ge 0}$ is a sequence in $L^1 (\langle v \rangle^\gamma)$ 
such that $g_n \to g$ and $\ron{L}(g_n) \to G$ in $L^1$, 
then $g \in L^1 (\langle v \rangle^\gamma)$ and $g_n \to g$ in $L^1 (\langle v \rangle^\gamma)$ 
by~\eqref{eq:compprec}. 
Then by~\eqref{eq:gain:poids}, \eqref{eq:dec:M}, \eqref{eq:nu} 
we deduce that $\ron{L}(g_n) \to \ron{L}(g)$ in $L^1$,  
which implies $G \equiv \ron{L}(g)$.
\end{proof}

\subsection{Estimates on $L$}
We recall here some classical properties of $L$.
 \begin{proposition}\label{prop:ferm:L}
  \begin{enumerate}
  \item[(i)] 
  The operators $L^+$ and $L^*$ 
  are bounded on $L^2(M)$. 
  \item[(ii)] The operator $L^\nu$, with domain $L^2 (\langle v \rangle^{2\gamma} M)$, 
  is closed on $L^2(M)$. 
  \item[(iii)] 
  The operator $L$, with domain $L^2 (\langle v \rangle^{2\gamma} M)$, 
  is closed on $L^2(M)$. 
  \end{enumerate} 
 \end{proposition}

\begin{proof}[Proof of proposition~\ref{prop:ferm:L}]
Point~(i) was proved in \eqref{eq:L^+bdd} (see also~\cite[Section 4]{Grad63} 
or~\cite[Chapter 7, Section 7.2]{CIP}). 
Point~(ii) is exactly similar to point~(iii) in Proposition~\ref{prop:ferm}. Point~(iii) is 
a consequence of the fact that a bounded perturbation of a closed operator is  
closed (see~\cite[Chapter~3, Section~5.2]{Kato}). 
\end{proof}

\section{Localization of the spectrum}
\label{spec}
\setcounter{equation}{0}

In this section we determine the spectrum of $\ron{L}$. 
As we do not have a hilbertian structure anymore, new technical difficulties 
arise with respect to the study of $L$. 
Nevertheless the localization of the essential spectrum is  
based on a similar argument as for $L$, 
namely the use of a variant of Weyl's Theorem for relatively compact perturbation. 
We assume in this section 
that the collision kernel $B$ satisfies \eqref{eq:prod}, \eqref{eq:hyprad}, \eqref{eq:grad}.

\subsection{Spectrum of $L$}

Before going into the study of the spectrum of $\ron{L}$ we 
state well-known properties on the spectrum of $L$. We recall that 
the discrete spectrum 
is defined as the set of eigenvalues isolated in 
the spectrum and with finite multiplicity, while the essential 
spectrum is defined as the complementary set in the spectrum of 
the discrete spectrum. 

 \begin{proposition}\label{prop:specL}
 The spectrum of $L$ is composed of an essential spectrum part,  
 which is $-\nu(\R^N) = (-\infty,-\nu_0]$, plus discrete eigenvalues 
 on $(-\nu_0,0]$, that can only accumulate at $-\nu_0$. 
 \end{proposition}
\begin{proof}[Proof of Proposition~\ref{prop:specL}]
The operator $L^c=L^+-L^*$ is compact on the Hilbert space $L^2(M)$ (see below). Thus 
Weyl's Theorem for self-adjoint operators (cf.~\cite[Chapter~4, Section~5]{Kato}) implies that 
 \[ \Sigma_e(L) = \Sigma_e (L^\nu) = (-\infty,-\nu_0]. \]
Since the operator is self-adjoint, the remaining part of the spectrum 
(that is the discrete spectrum) is included in 
$\R \cap (\C \setminus \Sigma_e(L)) = (-\nu_0,+\infty)$. 
Finally since the Dirichlet form is nonpositive, the discrete 
spectrum is also included in $\R_-$, which concludes the proof. 

Concerning the proof of the compactness of $L^c$, we shall briefly recall 
the arguments. The original proof is due to Grad \cite[Section 4]{Grad63} 
(in dimension $3$ for cutoff hard potentials). It was 
partly simplified in \cite[Chapter 7, Section 2, Theorem 7.2.4]{CIP} 
(in dimension $3$ for hard spheres). It relies on the Hilbert-Schmidt theory for integral 
operators (see \cite[Chapter 5, Section 2.4]{Kato}). We give here a version 
valid for cutoff hard potentials (under our assumptions \eqref{eq:prod}, 
\eqref{eq:hyprad}, \eqref{eq:grad}), in any dimension $N \ge 2$.

Let us first consider 
the compactness of $L^+$. First it was proved within Proposition \ref{prop:cvgL} the convergence 
 \[ \big\| {\bf 1}_{\{|\cdot|\le R\}} L^+ - L^+\big\|_{L^2(M)} \xrightarrow[]{R \to +\infty} 0. \]
Hence it is enough to prove the compactness of ${\bf 1}_{\{|\cdot|\le R\}} L^+$ 
for any $R>0$. Second if one defines 
 \[ L^+ _\var (h) := \int_{\R^N \times \ens{S}^{N-1}} {\bf 1}_{\{|v-v'|\ge \var\}} \, 
    \Phi(|v-v_*|) \, b(\cos \theta) \, 
 \left[ h' + h' _* \right] \, M_*  \, dv_* \, d\sigma, \]
the same computations of the kernel as in Proposition \ref{prop:cvgL} show that 
 \[ L^+ _\var (h) (v) = M^{-1/2} (v) \, \int_{u \in \R^N} k_\var(u,v) \, \left( h(u) M^{1/2}(u) \right) \, du \]
where $k_\var$ satisfies
  \[ k_\var(u,v) \le C \, {\bf 1}_{\{|u-v|\ge \var\}} \,  
   |u-v|^{1+\gamma-N} \, \exp \left[-\frac{|u-v|^2}{4}\right]. \]
By Young's inequality we deduce that 
  \[ \big\| L^+ _\var - L^+\big\|_{L^2(M)} \le C \, \left\| {\bf 1}_{\{|\cdot|\le \var\}} \,  
   |\cdot|^{1+\gamma-N} \, \exp \left[-\frac{|\cdot|^2}{4}\right] \right\|_{L^1}  
      \xrightarrow[]{\var \to 0} 0 \]
since the function $|\cdot|^{1+\gamma-N} \, \exp \left[-\frac{|\cdot|^2}{4}\right]$ 
is integrable at $0$. Hence it is enough to prove the compactness of 
${\bf 1}_{\{|\cdot|\le R\}} L^+ _\var$ for any $R,\var>0$. But as 
  \[ \int_{\R^N \times \R^N} {\bf 1}_{\{|v|\le R\}} \, {\bf 1}_{\{|u-v|\ge \var\}} \, 
              \left( |u-v|^{1+\gamma-N} \, \exp \left[-\frac{|u-v|^2}{4}\right] \right)^2 \, du \, dv < +\infty, \]
it is a Hilbert-Schmidt operator. This concludes the proof for $L^+$. For 
$L^*$, straightforward computations show that 
  \[ L^* (h) (v) = M^{-1/2} (v) \, \int_{u \in \R^N} k^*(u,v) \, \left( h(u) M^{1/2}(u) \right) \, du \]
with a kernel $k^*$ satisfying
  \[ k^*(u,v) \le C \, |u-v|^{\gamma} \, \exp \left[-\frac{|u|^2+|v|^2}{2}\right]. \]
This shows by inspection that $L^*$ is a Hilbert-Schmidt operator. 
\end{proof}

\subsection{Essential spectrum of $\ron{L}$}
 
Now let us turn to $\ron{L}$. We prove that 
the operator $\ron{L}^c$ is relatively compact with respect to $\ron{L}^\nu$. 
The main ingredients are the regularity 
estimates~\eqref{eq:reg:gain} and~\eqref{eq:reg:M},  
related to the ``almost convolution'' structure  
of the non-local term. We first deal with the approximate operator.  

 \begin{lemma}\label{lem:cpctdelta}
 For all $\delta \in (0,1)$, the 
 operator $\ron{L}^c _\delta$ is compact on $L^1$. 
 \end{lemma}
 
\begin{proof}[Proof of Lemma~\ref{lem:cpctdelta}]
We fix $\delta \in (0,1)$. We have to prove that for any sequence 
$(g_n)_{n \ge 0}$ bounded in $L^1$, the sequence 
$\left(\ron{L}^c _\delta (g_n)\right)_{n \ge 0}$ 
has a cluster point in $L^1$. 
The regularity estimates~\eqref{eq:reg:gain} 
and~\eqref{eq:reg:M} on $\ron{L}^+ _\delta$ and $\ron{L}^* _\delta$ 
imply that the sequence $\left(\ron{L}^c _\delta (g_n)\right)_{n \ge 0}$ 
is bounded in $W^{1,1}(\R^N)$. Then we can apply the Rellich-Kondrachov Theorem 
(see~\cite[Chapter~9, Section~3]{Brezis})  
on any open ball $B(0,K) \subset \R^N$ for $K \in \N^*$. 
It implies that the restriction of the sequence $\left(\ron{L}^c _\delta (g_n)\right)_{n \ge 0}$ 
to this ball is relatively compact in $L^1$. 
By a diagonal process with respect to the parameter $K \in \N^*$, 
we can thus extract a subsequence converging 
in $L^1 _{\mbox{{\scriptsize loc}}} (\R^N)$. The decay estimates~\eqref{eq:dec:gain} 
and~\eqref{eq:dec:M} then ensure a tightness control 
(uniform with respect to $n$), which implies that the 
convergence holds in $L^1$. This ends the proof.  
\end{proof}

Then by closeness of the relative compactness property, we deduce for $\ron{L}^c$
 \begin{lemma}\label{lem:cpct}
 The operator $\ron{L}^c$ is relatively compact with respect to $\ron{L}^\nu$. 
 \end{lemma}

\begin{proof}[Proof of Lemma~\ref{lem:cpct}]
Thanks to the estimate~\eqref{eq:nu}, it is equivalent to prove 
that for any sequence $(g_n)_{n \ge 0}$ bounded in 
$L^1 (\langle v \rangle^\gamma)$, the sequence $\left(\ron{L}^c (g_n)\right)_{n \ge 0}$ 
has a cluster point in $L^1$. 
As $L^1$ is a Banach space it is 
enough to prove that a subsequence of 
$\left(\ron{L}^c (g_n)\right)_{n \ge 0}$ has the Cauchy property. 
Let us choose a sequence $\delta_k \in (0,1)$ 
decreasing to $0$. Thanks to the 
previous lemma and a diagonal process, 
we can find an extraction $\varphi$ such that for all 
$k \ge 0$ the sequence 
$\left(\ron{L}^c _{\delta_k} (g_{\varphi(n)})\right)_{n \ge 0}$ 
converges in $L^1$. 
Then for a given $\var >0$, we first choose $k \in \N$ such that 
 \[ \forall \, n \ge 0, \ \ \ \| \ron{L}^c _{\delta_k} (g_n) 
       - \ron{L}^c (g_n) \|_{L^1} \le \var/4 \] 
which is possible thanks to Proposition~\ref{prop:cvgronL} 
and the uniform bound on the $L^1 (\langle v \rangle^\gamma)$ norm of the sequence $(g_n)_{n \ge 0}$. 
Then we choose $n_{k,\var}$ such that 
 \[  \forall \, m,n \ge n_{k,\var}, \ \ \ 
     \| \ron{L}^c _{\delta_k} (g_{\varphi(m)}) 
     - \ron{L}^c _{\delta_k} (g_{\varphi(n)}) \| \le \var/2 \] 
since the sequence 
$\left(\ron{L}^c _{\delta_k} (g_{\varphi(n)})\right)_{n \ge 0}$ 
converges in $L^1$. Then by triangular inequality, we get 
 \[  \forall \, m,n \ge n_{k,\var}, \ \ \ 
       \|\ron{L}^c (g_{\varphi(m)}) - 
       \ron{L}^c (g_{\varphi(n)}) \| \le \var, \] 
which concludes the proof.
\end{proof}

The next step is the use of a variant of Weyl's Theorem. 

 \begin{proposition}\label{prop:ess}
 The essential spectrum of the operator $\ron{L}$ 
 is $-\nu(\R^N)=(-\infty,-\nu_0]$.
 \end{proposition}
 
\begin{proof}[Proof of Proposition~\ref{prop:ess}]
We shall use here the classification of the spectrum by the 
Fredholm theory. Indeed in the case of non hilbertian operators, Weyl's  
Theorem does not imply directly the stability of the essential 
spectrum under relatively compact perturbation, but only the 
stability of a smaller set, namely the complementary in the 
spectrum of the Fredholm set (see below). 
We refer for the objects and results   
to~\cite[Chapter~4, Section~5.6]{Kato}. 

Given an operator $T$ on a Banach space 
${\cal B}$ and a complex number $\xi$, 
we define $\mbox{nul}(\xi)$ as the dimension 
of the null space of $T-\xi$, 
and $\mbox{def}(\xi)$ as the codimension 
of the range of $T- \xi$. 
These numbers belong to $\N \cup \{+\infty\}$. 
A complex number $\xi$ belongs to the 
resolvent set if and only if 
$\mbox{nul}(\xi) = \mbox{def}(\xi) =0$. 
Let $\Delta_F (T)$ be the set of all complex 
numbers such that $T-\xi$ is Fredholm 
({\it i.e.} $\mbox{nul}(\xi) < +\infty$ and 
$\mbox{def}(\xi) < +\infty$). 
This set includes the resolvent set. 
Let $E_F (T)$ be the complementary set 
of $\Delta_F (T)$ in $\C$, in short the set of complex numbers 
$\xi$ such that $T-\xi$ is not Fredholm. 
From~\cite[Chapter 4, Section~5.6, Theorem~5.35 and footnote]{Kato}, 
the set $E_F$ is preserved under relatively compact perturbation. 

We apply this result to the perturbation 
of $-\ron{L}^\nu$ by $\ron{L}^c$, which is relatively compact by Lemma~\ref{lem:cpct}. 
As $E_F (\ron{L}^\nu) = -\nu(\R^N) = (-\infty, \nu_0]$, 
we deduce that $E_F (\ron{L})= (-\infty,-\nu_0]$ and so 
$\Delta_F (\ron{L}) = \ens{C} \setminus (-\infty,-\nu_0]$. 

Thus it remains to prove that the Fredholm set of $\ron{L}$ 
contains only the discrete spectrum plus the resolvent set. 
By~\cite[Chapter~4, Section~5.6]{Kato}, 
the Fredholm set $\Delta_F$ is the union of a countable number 
of components $\Delta_n$ (connected open sets) 
on which $\mbox{nul}(\xi)$ and $\mbox{def}(\xi)$ 
are constant, except for a (countable) set of isolated values of $\xi$. 
Moreover the boundary $\partial \Delta_F$ of the 
set $\Delta_F$ as well as the boundaries $\partial \Delta_n$ 
of the components $\Delta_n$ all belong to the set $E_F$. 
As in our case the Fredholm set 
$\Delta_F (\ron{L})= \C \setminus (-\infty,-\nu_0]$ 
is connected, it has only one component. 
It means that $\mbox{nul}(\xi)$ and $\mbox{def}(\xi)$ 
are constant on $\C \setminus (-\infty,-\nu_0]$, 
except for a (countable) set of isolated values of $\xi$. 

Let us prove now that these constant values are 
$\mbox{nul}(\xi) = \mbox{def}(\xi) =0$. It 
will imply the result, since a complex number $\xi$ such that 
$\mbox{nul}(\xi) = \mbox{def}(\xi) =0$ belongs to the resolvent set, 
and a complex number $\xi$, isolated in the spectrum, that belongs to 
the Fredholm set, satisfies $\mbox{nul}(\xi) < +\infty$ and 
$\mbox{def}(\xi) < +\infty$, and is exactly 
a discrete eigenvalue with finite multiplicity. 

As the numbers $\mbox{nul}(\xi)$ and $\mbox{def}(\xi)$ are constant in 
$\Delta_F (\ron{L})= \C \setminus \left( (-\infty,\nu_0] \cup {\cal V} \right)$ 
(${\cal V}$ denotes a (countable) set of isolated complex numbers), it is enough to prove that 
there is an uncountable set of complex numbers in 
$\C \setminus (-\infty,\nu_0]$ such that 
$\mbox{nul}(\xi) = \mbox{def}(\xi) =0$.

By using \eqref{eq:auxcomp} and \eqref{eq:dec:M} we have 
  \begin{equation}\label{eq:rel}
  \forall \, g \in L^1 (\langle v \rangle^\gamma), \ \ \ 
  \| \ron{L}^c (g) \|_{L^1} \le 
     C \, \| g\|_{L^1} + \frac{n_0}2 \, \| g\|_{L^1 (\langle v \rangle^\gamma)} 
  \le C \, \| g\|_{L^1} + \frac12 \, \|\nu g\|_{L^1}
  \end{equation}
for some explicit constant $C$. 
Now we choose $r_0 \in \R_+$ big enough such that 
 \[ \forall \, r \ge r_0, \ \ \ \frac{C}{\nu_0 + r} + \frac12 < 1. \] 
The multiplication operator $-(\nu+r)$ is bijective from 
$L^1 (\langle v \rangle^\gamma)$ to $L^1$ (since $\nu+r > \nu_0>0$). 
The inverse linear operator 
is the multiplication operator $-(\nu+r)^{-1}$, it is defined 
on $L^1$ and bounded by $\|-(\nu+r)^{-1}\|_{L^\infty} = (\nu_0 +r) ^{-1}$. 
Its range is $L^1 (\langle v \rangle^\gamma)$. 
Hence the linear operator $\ron{L}^c ( -(\nu+r)^{-1} \cdot)$ is 
well-defined, and, thanks to~\eqref{eq:rel} it is bounded with 
a norm controlled by $C/(\nu_0 + r) + 1/2$,  
which is strictly less than $1$ for $r \ge r_0$. 
Thus for $r \ge r_0$, the operator $\mbox{Id}+\ron{L}^c ( -(\nu+r)^{-1} \cdot)$ 
is invertible with bounded inverse, and as the operator $ -(\nu+r) \cdot$ 
is also invertible with bounded inverse for $r \ge 0$, by composition we deduce that  
  \[ \Big( \mbox{Id}+\ron{L}_c ( -(\nu+r)^{-1} \cdot) \Big) 
     \circ \left( -(\nu+r) \cdot \right) 
     = \ron{L}_c -(\nu+r) \] 
is invertible with bounded inverse. It means that $[r_0,+\infty)$ belongs 
to the resolvent set, {\it i.e.}  
$\mbox{nul}(r) = \mbox{def}(r) =0$ for all $r \ge r_0$, 
which concludes the proof. 
\end{proof}

\subsection{Discrete spectrum of $\ron{L}$}

In order to localize the discrete eigenvalues, 
we will prove that the eigenvectors associated with these eigenvalues 
decay fast enough at infinity to be in fact multiple of the 
eigenvectors of $L$. This implies that 
these eigenvalues belong to the discrete spectrum of $L$  
and gives new geometrical informations on them: they lie in the intervalle $(-\nu_0,0]$ 
with the only possible cluster point being $-\nu_0$.  
Moreover, explicit estimates on the spectral gap of $\ron{L}$ follow by~\cite{BaMo}. 

  \begin{proposition}\label{prop:decvp}
  The operators $\ron{L}$ and $L$ 
  have the same discrete eigenvalues with the same 
  multiplicities. Moreover the eigenvectors of $\ron{L}$ associated  
  with these eigenvalues are given by those of $L$ associated 
  with the same eigenvalues, multiplied by $m^{-1} M$. 
  \end{proposition}

\Remarks 
This result implies in particular that the (finite dimensional) 
algebraic eigenspaces of the discrete eigenvalues of $\ron{L}$ 
do not contain any Jordan block (i.e. their algebraic 
multiplicity equals their geometric multiplicity, see the 
definitions in \cite[Chapter 3, Section 6]{Kato}) 
as it is the case for the self-adjoint operator $L$. 
\medskip


\begin{proof}[Proof of Proposition~\ref{prop:decvp}]
Let us pick $\lambda$ a discrete eigenvalue of $\ron{L}$. 
The associated eigenspace has finite dimension since 
the eigenvalue is discrete. Let us consider 
a Jordan block of $\ron{L}$ on this eigenspace, spanned 
in the canonical form by the basis $(g_1,g_2,\dots,g_n)$. It means that 
 \[ \ron{L} (g_1) = \lambda \, g_1 \]
and for all $2 \le i \le n$, 
 \[ \ron{L} (g_i) = \lambda \, g_i + g_{i-1}. \]
  
As $\lambda$ does not belong to the essential 
spectrum of $\ron{L}$, we know from~Proposition~\ref{prop:ess}  
that $\lambda \not \in (-\infty,-\nu_0]$. 
Let us call $d_\lambda >0$ the distance 
between $\lambda$ and $(-\infty,-\nu_0]$. It is straightforward that 
  \[ \forall \, v \in \R^N, \ \ \ 
      |\nu (v) + \lambda| \ge d_\lambda \]
and by~\eqref{eq:nu} there is $d_\lambda '>0$ such that
  \[ \forall \, v \in \R^N, \ \ \ 
      |\nu (v) + \lambda| \ge d_\lambda ' \, \langle v \rangle^\gamma. \]

Let us prove by finite induction that for all $1 \le i \le n$ 
we have $m M^{-1} g_i \in L^2 (\langle v \rangle^{2\gamma} M)$.
We write 
  \[ \ron{L} -\lambda = \left( \ron{L}^+ _\delta - \ron{L}^* \right) 
                      - \left( \nu + \lambda - \ron{L}^+ + \ron{L}^+ _\delta \right) 
                       =: A_\delta - B_\delta.  \]
Both part $A_\delta$ and $B_\delta$ of this decomposition are well-defined on 
$L^1(\langle v \rangle^\gamma)$. Moreover we shall prove that when 
$\delta$ is small enough, $B_\delta$ is bijective from $L^1(\langle v \rangle^\gamma)$ 
to $L^1$ with bounded inverse, and also 
that its restriction $(B_\delta)_|$ to $L^2(\langle v \rangle^{2\gamma} m^2 M^{-1}) \subset L^1$ is   
bijective from $L^2(\langle v \rangle^{2\gamma} m^2 M^{-1})$ to $L^2(m^2 M^{-1})$ 
with bounded inverse. 

We pick $\delta > 0$ such that 
  \begin{equation}\label{choixdelta}
  \forall \, v \in \R^N, \ \ \ C_1 (\delta) \le \frac{d_\lambda '}{2} 
                 \ \ \mbox{ and } \ \  C_2(\delta) \le  \frac{d_\lambda}{2} 
  \end{equation}
where $C_1(\delta)$ and $C_2(\delta)$ are defined in Propositions~\ref{prop:cvgronL} 
and~\ref{prop:cvgL}. Then we write 
 \[ B_\delta = \Big( \mbox{Id} - (\ron{L}^+ - \ron{L}^+ _\delta) ( (\nu+\lambda)^{-1} \cdot) \Big) 
    \circ \left( (\nu+\lambda) \cdot \right). \]
As $(\ron{L}^+ - \ron{L}^+ _\delta) ( (\nu+\lambda)^{-1} \cdot)$ is bounded in $L^1$ with norm 
less than $1/2$ thanks to~\eqref{choixdelta}, we have that 
$\mbox{Id} - (\ron{L}^+ - \ron{L}^+ _\delta) ( (\nu+\lambda)^{-1} \cdot)$ 
is bijective from $L^1$ to $L^1$ (with bounded inverse). 
As $(\nu+\lambda) \cdot$ is bijective from $L^1(\langle v \rangle^\gamma)$ to $L^1$ 
(with bounded inverse), we deduce that $B_\delta$ is bijective from 
$L^1(\langle v \rangle^\gamma)$ to $L^1$ (with bounded inverse). 

Then we remark that 
  \[ \big\| (\ron{L}^+ - \ron{L}^+ _\delta) \big\|_{L^2(m^2 M^{-1})} = \big\|L^+ - L^+ _\delta\big\|_{L^2(M)} \]
thanks to the formula \eqref{eq:auxdefronL^+}, \eqref{eq:auxdefL^+}, 
\eqref{eq:auxdefronL^+delta}, \eqref{eq:auxdefL^+delta} 
for $\ron{L}^+$, $L^+$, $\ron{L}^+ _\delta$ and $L^+ _\delta$. 
Hence $(\ron{L}^+ - \ron{L}^+ _\delta) ( (\nu+\lambda)^{-1} \cdot)$ 
is bounded in $L^2(m^2M^{-1})$ with norm less than $1/2$ thanks to~\eqref{choixdelta}, 
and we deduce that $\left( \mbox{Id} - (\ron{L}^+ - \ron{L}^+ _\delta) ( (\nu+\lambda)^{-1} \cdot) \right)$ 
is bijective from $L^2(m^2M^{-1})$ to $L^2(m^2M^{-1})$ (with bounded inverse). 
As the multiplication operator $\nu +\lambda$ is bijective from $L^2(\langle v \rangle^{2\gamma} m^2M^{-1})$ 
to $L^2(m^2M^{-1})$ (with bounded inverse), 
we deduce that $(B_\delta)_|$ is bijective from $L^2(\langle v \rangle^{2\gamma} m^2M^{-1})$ 
to $L^2(m^2M^{-1})$ (with bounded inverse). 
 
For the initialization, we write the eigenvalue equation on $g_1$ in the form 
  \begin{equation}\label{eqvp} 
  B_\delta (g_1) = A_\delta(g_1). 
  \end{equation}
Thanks to the decay estimates~\eqref{eq:dec:gain} and~\eqref{eq:dec:M}, $A_\delta(g_1)$ 
belongs to $L^2(m^2M^{-1}) \subset L^1$, and thus it implies that the unique 
pre-image of $A_\delta(g_1)$ by $B_\delta$ in $L^1$ belongs to 
$L^2(\langle v \rangle^{2\gamma} m^2M^{-1})$, thanks to the invertibilities 
of $B_\delta$ and $(B_\delta)_|$ proved above. 
Hence $g_1 \in L^2(\langle v \rangle^{2\gamma} m^2M^{-1})$.  

Now let us consider the other vectors of the Jordan block:  
we pick $2 \le i \le n$ and we  
suppose the result to be true for $g_{i-1}$. 
Then $g_i$ satisfies 
 \[ B_\delta (g_i) = A_\delta(g_i) - g_{i-1} \]
and with the same argument as above together with the fact that 
$g_{i-1} \in L^2(\langle v \rangle^{2\gamma} m^2M^{-1})$, one concludes straightforwardly. 

As a consequence, for any $1 \le i \le n$, $g_i$ belongs to 
$L^2(\langle v \rangle^{2\gamma} m^2 M^{-1})$ and thus $m M^{-1} g_i$ 
belong to the space $L^2(\langle v \rangle^{2\gamma} M)$, {\it i.e.} the domain of $L$.  
Hence $\lambda$ is necessarily an eigenvalue of $L$, and 
the eigenspace associated with $\lambda$ of the operator $\ron{L}$ is included in 
the one of $L$ multiplied by $m^{-1} M$. As the converse 
inclusion is trivially true, this ends the proof. 
\end{proof}

\begin{figure}[h]
  \epsfysize=4cm
  $$\epsfbox{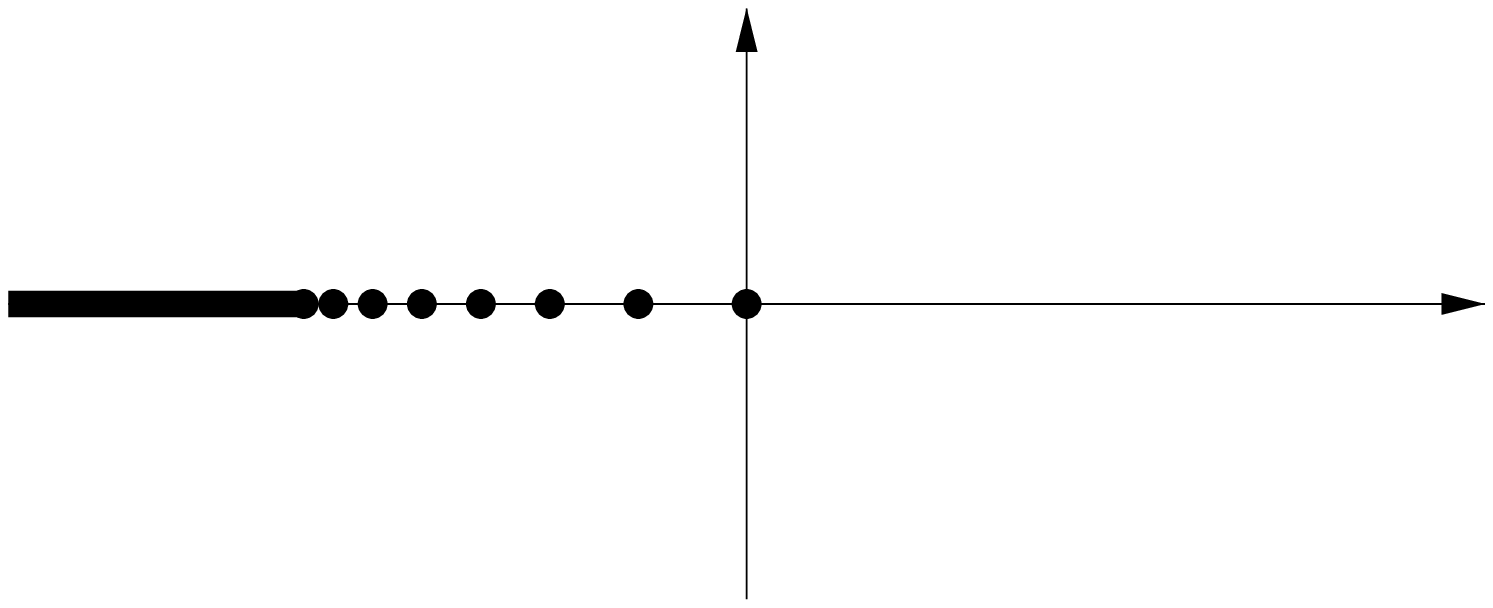}$$
  \caption{Spectrum of $\ron{L}$ in $L^1$}\label{fig:spec:ronL}
  \end{figure}

To conclude this section, we give in Figure~\ref{fig:spec:ronL} 
the complete picture of the spectrum of $\ron{L}$ in $L^1$, 
which is the same as the spectrum of $L$ in $L^2(M)$ 
(using Proposition~\ref{prop:ess} and Proposition~\ref{prop:decvp}).

\section{Trend to equilibrium}
\label{cvg}
\setcounter{equation}{0}

This section is devoted to the proof of the main Theorem~\ref{theo:cvg}.  
We consider $f$ a solution (in $L^1(\langle v \rangle^2)$) of the Boltzmann equation~\eqref{eq:base}. 
The equation satisfied by the perturbation of equilibrium $g=m^{-1} \, (f-M)$ is
  \[ \frac{\partial g}{\partial t} = \ron{L}(g) + \Gamma(g,g) \]
with 
  \[ \Gamma(g,g) = m^{-1} \, Q(m g, m g). \] 
The null space of $\ron{L}$ is given by the one of $L$ multiplied 
by $m^{-1} \, M$ (cf. Proposition~\ref{prop:decvp}), {\it i.e.} 
the following $(N+2)$-dimensional vector space:
  \[ N(\ron{L}) = 
    \mbox{Span} \left\{ m^{-1} \, M, m^{-1} \, M \, v_1, \dots, m^{-1} \, M \, v_N , m^{-1} \, M \, |v|^2 \right\}. \] 

Let us introduce the following complementary set of 
$N(\ron{L})$ in $L^1$:
 \[ \cal{S} = \left\{ g \in L^1, \hspace{0.3cm}  
     \int_{\R^N} m\, g \, \phi \, dv = 0, \ \ 
     \phi(v) = 1, v_1, \dots, v_N, |v|^2 \right\}.  \]   
Since
 \[ \int_{\R^N} \ron{L}(g) \, m \, \phi \, dv = 
    \int_{\R^N \times \R^N \times \ens{S}^{N-1}} \Phi \, b \, m \, g \, M_* \, 
           \left[ \phi' + \phi' _* - \phi - \phi_* \right] \, dv \, dv_* \, d\sigma, \]
 \[ \int_{\R^N} \Gamma(g,g) \, m \, \phi \, dv = 
    \int_{\R^N \times \R^N \times \ens{S}^{N-1}} \Phi \, b \, m \, g \, m_* \, g_* \, 
           \left[ \phi' + \phi' _* - \phi - \phi_* \right] \, dv \, dv_* \, d\sigma, \]
we see that
  \[ \ron{L} (L^1(\langle v \rangle^\gamma)) \subset \cal{S}, \hspace{1cm} 
    \Gamma (L^1(\langle v \rangle^\gamma),L^1(\langle v \rangle^\gamma)) \subset \cal{S}. \] 
As $g_0 \in \cal{S}$ since $f$ and $M$ have the same mass, momentum and energy, 
we can restrict the evolution equation to $\cal{S}$ and 
thus we shall consider in the sequel the operator $\ron{L}$ 
restricted on $\cal{S}$, which we denote by 
$\ron{\tilde{L}}$. The spectrum of $\ron{\tilde{L}}$ is 
given by the one of $\ron{L}$ minus the $0$ eigenvalue. 

Similarly we define in $L^2(M)$ the following stable 
complementary set of the kernel $N(L)$ of $L$ ($N(L)$ was defined in~\eqref{noyauL}) 
  \[ S = \left\{ h \in L^2(M), \hspace{0.5cm}  
      \int_{\R^N} h \, \phi \, M \, dv = 0, \ \ 
      \phi(v) = 1, v_1, \dots, v_N, |v|^2 \right\} \] 
which is formally related to $\cal{S}$ by $S= m \, M^{-1} \, \cal{S}$.
We define the restriction $\tilde{L}$ of $L$ on the stable set $S$, 
whose spectrum is given by the one of $L$ minus the $0$ eigenvalue. 


\subsection{Decay estimates on the evolution semi-group} 
\label{subsec:decsemi}

Compared to the classical strategy to obtain decay estimates on the evolution semi-group of 
$L$, here the estimate on the Dirichlet form will be replaced by an estimate on 
the norm of the resolvent and the self-adjointness property will replaced by the 
sectorial property. 

We denote $\Sigma = \Sigma(L) = \Sigma(\ron{L})$ and 
$\Sigma_e = \Sigma_e (L) = \Sigma_e (\ron{L})$. 
For $\xi \not\in \Sigma$, we denote  
$\ron{R}(\xi) = ( \ron{L} -\xi )^{-1}$ the resolvent of 
$\ron{L}$ at $\xi \in \C$, and 
$R(\xi) = (L -\xi )^{-1}$ the resolvent of 
$L$ at $\xi \in \C$. 
$\ron{R}(\xi)$ is a bounded operator on $L^1$, 
$R(\xi)$ is a bounded operator on $L^2(M)$. 
We have the following estimate on the norm of $\ron{R}(\xi)$:

  \begin{proposition}\label{prop:normR}
  There are explicit constants $C_8$, $C_9>0$ depending 
  only on the collision kernel and on a lower bound on 
  $\mbox{{\em dist}}(\xi,\Sigma_e)$ such that 
    \begin{equation}\label{eq:normR}
    \forall \, \xi \not \in \Sigma, \ \ \ 
    \| \ron{R}(\xi) \|_{L^1} \le C_8 + C_9 \, \| R(\xi) \|_{L^2(M)}.
    \end{equation} 
  \end{proposition}

\begin{proof}[Proof of Proposition~\ref{prop:normR}]
%
Let us introduce a right inverse of the operator $(\ron{L}-\xi)$: 
let us pick $\delta > 0$ such that 
  \begin{equation}\label{choixdeltabis}
  \forall \, v \in \R^N, \ \ \ C_1 (\delta) \le \frac{\nu(v)+\xi}{2 \, \langle v \rangle^\gamma} 
  \end{equation}
(note that this choice only depends on the collision kernel and a lower bound on $\mbox{dist}(\xi,\Sigma_e)$). 
Then we shall use the same argument as in the proof of Proposition~\ref{prop:decvp} 
to prove that the operator 
  \[ B_\delta(\xi) = \ron{L}^\nu + \xi - \left( \ron{L}^+ - \ron{L}^+ _\delta \right) \]
is bijective from $L^1(\langle v \rangle^\gamma)$ to $L^1$ and its inverse has its norm bounded by 
  \[ \|B_\delta(\xi) ^{-1} \|_{L^1} \le \frac{2}{\mbox{dist}(\xi,\Sigma_e)}. \]
Indeed once the invertibility is known, the bound on the inverse is given by 
  \begin{eqnarray*}
  \forall \, g \in L^1 (\langle v \rangle^\gamma), \ \ \ \|B_\delta(\xi)(g)\|_{L^1} 
        &\ge& \|(\nu+\xi)g\|_{L^1} - \|(\ron{L}^+ - \ron{L}^+ _\delta)(g)\|_{L^1} \\
        &\ge& 
             \frac{\mbox{dist}(\xi,\Sigma_e)}{2} \, \|g\|_{L^1}
  \end{eqnarray*}
where we have used Proposition \ref{prop:cvgronL} and the bound \eqref{choixdeltabis} 
on $C_1 (\delta)$. To prove the invertibility, we write  
  \[ B_\delta(\xi) = \Big( \mbox{Id} +(\ron{L}^+ - \ron{L}^+ _\delta) ( -(\nu+\xi)^{-1} \cdot) \Big) 
    \circ \left( -(\nu+\xi) \cdot \right). \] 
Since $(\ron{L}^+ - \ron{L}^+ _\delta) ( -(\nu+\xi)^{-1} \cdot)$ is well-defined and 
bounded on $L^1$ with norm less than $1/2$ thanks to~\eqref{choixdeltabis}, 
we have that $\mbox{Id} +(\ron{L}^+ - \ron{L}^+ _\delta) ( -(\nu+\xi)^{-1} \cdot)$ 
is bijective from $L^1$ to $L^1$ (with bounded inverse), and as 
$-(\nu+\xi) \cdot$ is bijective from $L^1(\langle v \rangle^\gamma)$ to $L^1$ (with bounded inverse), 
we deduce the result. 

Now we denote
  \[ A_\delta = \ron{L}^+ _\delta - \ron{L}^*. \]
This operator is bounded on $L^1$ and satisfies, thanks to the 
estimates~\eqref{eq:dec:gain} and~\eqref{eq:dec:M},
  \[ \forall \, v \in \R^N, \ \ \ |A_\delta(g)(v)| \le C \, \|g\|_{L^1} \, M^\theta \]
for any $\theta \in [0,1)$ and some explicit constant $C$ depending 
on the choice of $\delta$ above, on the collision kernel, and on $\theta$. 

The operator $\ron{L}-\xi$ writes 
  \[ \ron{L}-\xi = A_\delta - B_\delta(\xi) \]
and we define the following operator 
  \[ I(\xi) = - B_\delta(\xi) ^{-1} 
       + (m^{-1} M) \, R(\xi) \left[ (m M^{-1}) \, A_\delta \, B_\delta(\xi) ^{-1} \right] \]
(note that here $R$ is the resolvent of $L$). 
Let us first check that this operator is well-defined and bounded on $L^1$: 
for $g \in L^1$, we have $B_\delta(\xi) ^{-1}(g) \in L^1$ and as (choosing $\theta=3/4$ for instance)
 \[ \big| A_\delta \big( B_\delta(\xi) ^{-1}(g)\big) (v) \big| 
    \le C \, \big\|B_\delta(\xi) ^{-1}(g)\big\|_{L^1} \, M^{3/4}, \]
we have 
 \[ \left\| (m M^{-1}) \, A_\delta \big( B_\delta(\xi) ^{-1}(g)\big) \right\|_{L^2(M)} ^2 
    \le C^2 \, \big\|B_\delta(\xi) ^{-1}(g)\big\|_{L^1} ^2 \, \left( \int_{\R^N} m^2 \, M^{1/2} \, dv \right). \]
Thus 
  \[ R(\xi) \left[ (m M^{-1}) \, A_\delta \, B_\delta(\xi) ^{-1}(g) \right] \]
is well-defined and belongs to $L^2(M)$. Since by Cauchy-Schwarz
  \[ \big\|m^{-1} M \, h\big\|_{L^1} \le \big\|m^{-1} M^{1/2} \big\|_{L^2} \, \|h\|_{L^2(M)} \]
we deduce finally that $I(\xi)(g)$ is well-defined and belongs to $L^1$. 
Moreover, from the computations above we deduce 
 \[ \|I(\xi)(g)\|_{L^1} \le \|B_\delta(\xi) ^{-1}\| \, \Big( 1 + C \, \| R(\xi) \|_{L^2(M)} \,  
                \left\|m^{-1} M^{1/2} \right\|_{L^2} 
                               \, \left\|m M^{1/4} \right\|_{L^2} \Big) \, \|g\|_{L^1}. \]

Now let us check that $I(\xi)$ is a right inverse of $(\ron{L}-\xi)$: 
  \begin{multline*}
  (\ron{L}-\xi) \circ I(\xi) (g) = 
     (\ron{L}-\xi) \circ 
     \Big( - B_\delta(\xi) ^{-1} 
         + (m^{-1} M) \, R(\xi) \left[ (m M^{-1}) \, A_\delta \, B_\delta(\xi) ^{-1} \right]\Big) (g)   \\
  = \big(-A_\delta + B_\delta(\xi)\big) \circ B_\delta(\xi) ^{-1} (g)
    + (\ron{L}-\xi) \circ 
    \Big( (m^{-1} M) \, R(\xi) \left[ (m M^{-1}) \, A_\delta \, B_\delta(\xi) ^{-1}\right] \Big) (g)   \\
  = g - A_\delta \, B_\delta(\xi) ^{-1} (g) 
    + (\ron{L}-\xi) \circ 
    \Big( (m^{-1} M) \, R(\xi) \left[ (m M^{-1}) \, A_\delta \, B_\delta(\xi) ^{-1} \right] \Big) (g).
  \end{multline*}
Now as 
  \[  (m^{-1} M) \, R(\xi) 
             \left[ (m M^{-1}) \, A_\delta \, B_\delta(\xi) ^{-1} \right] (g)  \in L^2(m^2 \, M^{-1}) \]
we deduce that 
 \begin{multline*}
 (\ron{L}-\xi) \circ 
   \Big( (m^{-1} M) \, R(\xi) \left[ (m M^{-1}) \, A_\delta \, B_\delta(\xi) ^{-1} \right] \Big) (g) \\
   = m^{-1} \, M \, (L-\xi) \circ 
   R(\xi) \left[ (m M^{-1}) \, A_\delta \, B_\delta(\xi) ^{-1} \right] (g)  \\
   = m^{-1} \, M \, \left[ (m M^{-1}) \, A_\delta \, B_\delta(\xi) ^{-1} (g) \right] 
   = A_\delta \, B_\delta(\xi) ^{-1} (g).
 \end{multline*}
Collecting every term we deduce 
 \[ \forall \, g \in L^1, \quad (\ron{L}-\xi) \circ I(\xi) (g) = g. \]

Let us conclude the proof: whenever $\xi \not \in \Sigma$, $(\ron{L}-\xi)$ is bijective 
from $L^1(\langle v \rangle^\gamma)$ to $L^1$ with bounded inverse $\ron{R}(\xi)$, and 
we deduce that $\ron{R}(\xi)=I(\xi)$ and thus 
  \[ \|\ron{R}(\xi)\|_{L^1} \le \|B_\delta(\xi) ^{-1}\| \, \Big[ 1 + C \, \|R(\xi)\|_{L^2(M)} \,  
                 \left\|m^{-1} M^{1/2} \right\|_{L^2} 
                                \, \left\|m M^{1/4} \right\|_{L^2} \Big]. \]
As we have 
  \[ \|B_\delta(\xi) ^{-1}\| \le \frac{2}{\mbox{dist}(\xi,\Sigma_e)} \]
and the choice of $\delta$ (determining the constant $C$) depends only on the 
collision kernel and a lower bound on $\mbox{dist}(\xi,\Sigma_e)$, this ends the proof. 
\end{proof} 

Now we use this estimate in order to obtain some decay estimate on the 
evolution semi-group. We recall that $\lambda \in (0,\nu_0)$ denotes the 
spectral gap of $\ron{L}$ and $L$. 
  \begin{theorem}\label{theo:sg}
  The evolution semi-group of the operator $\ron{\tilde{L}}$ is well-defined on $L^1$, 
  and for any $0< \mu \le \lambda$, it satisfies the decay estimate
    \begin{equation}\label{eq:decest}
    \forall \, t \ge 0, \quad \big\| e^{t \ron{\tilde{L}}} \big\|_{L^1} \le C_{10} \, e^{-\mu t} 
    \end{equation}
  for some explicit constant $C_{10} >0$ depending only on the 
  collision kernel, on $\mu$, and a lower bound on $\nu_0 -\mu$. 
  \end{theorem}
\begin{proof}[Proof of Theorem~\ref{theo:sg}]
Let us pick $\mu \in (0,\lambda]$. We define in the complex plane the set 
  \[ {\cal A}_\mu = \left\{ \xi \in \C, \ \ \arg (\xi - \mu) \in 
                     \left[-\frac{3\pi}{4},\frac{3\pi}{4} \right] \ \ \mbox{ and } \ \  
                     \mbox{Re}(\xi) \le - \frac{\mu}{2} \right\}. \]
We shall prove the following lemma
  
  \begin{lemma}\label{lem:resolvent}
  There are explicit constants  
  $a,b>0$ depending on the collision kernel, on $\mu$,    
  and on a lower bound on $\nu_0 - \mu$, such that 
    \[ \forall \, \xi \in {\cal A}_\mu, \ \ \ \|{\cal R}(\xi)\|_{L^1} \le a + \frac{b}{|\xi -\mu|}. \]
  \end{lemma}
\begin{proof}[Proof of Lemma~\ref{lem:resolvent}]
We shall use Proposition~\ref{prop:normR}. In the Hilbert space $L^2(M)$, 
$L$ is self-adjoint and thus we have (see~\cite[Chapter~5, Section~3.5]{Kato}) 
  \[ \| R(\xi) \|_{L^2(M)} = \frac{1}{\mbox{dist}(\xi,\Sigma)}. \]
Hence Proposition~\ref{prop:normR} yields
  \[ \forall \, \xi \in {\cal A}_\mu, \quad  
                  \|{\cal R}(\xi)\|_{L^1} \le C_8 + \frac{C_9}{\mbox{dist}(\xi,\Sigma)} \]
with $C_8$ and $C_9$ depending on a lower bound on $\mbox{dist}(\xi,\Sigma_e)$. 
Then in the set ${\cal A}_\mu$, the lower bound on $\mbox{dist}(\xi,\Sigma_e)$ 
is straightforwardly controlled by a lower bound on $\nu_0 - \mu$, and we have immediately 
  \[ \mbox{dist}(\xi,\Sigma \setminus \{0\}) \ge \frac{|\xi - \mu|}{\sqrt{2}} \]
and $\mbox{dist}(\xi, \{0\}) =|\xi|$. Since for $\xi \in {\cal A}_\mu$, 
we have $|\xi-\mu|\le |\xi|$, we deduce that 
  \[  \forall \, \xi \in {\cal A}_\mu, \ \ \ \|{\cal R}(\xi)\|_{L^1} \le 
        C_8 + C_9 \, \max \left\{ \frac{\sqrt{2}}{|\xi - \mu|} , \frac{1}{|\xi|} \right\} 
                 \le a + \frac{b}{|\xi - \mu|}, \]
which concludes the proof.
\end{proof}
Now let us conclude the proof of the theorem. Let $t>0$ and $\eta \in (0,\pi/4)$. 
Let us consider $\Gamma$ a curve running, within ${\cal A}_\mu$, 
from infinity with $\arg(\xi)=\pi/2+\eta$ to infinity 
with $\arg(\xi) = -\pi/2 - \eta$, and the complex integral 
  \[ \frac{-1}{2 \pi i} \, \int_{\Gamma} e^{t \xi} \, {\cal R}(\xi) \, d\xi. \]
Thanks to the bound of Lemma~\ref{lem:resolvent}, the integral is absolutely convergent. As the 
curve encloses the spectrum of $\ron{L}$ minus $0$, {\it i.e.} the spectrum of $\ron{\tilde{L}}$,  
classical results from spectral analysis (see~\cite[Chapter~1, Section~3]{Henry} 
and~\cite[Chapter~9, Section~1.6]{Kato}) show that this integral defines the 
evolution semi-group $e^{t \ron{\tilde{L}}}$ of $\ron{\tilde{L}}$. 
Now we apply a classical strategy to obtain a decay estimate on the semi-group: 
we perform the change of variable $\xi = z/t - \mu$. Then $z$ describes a new path 
$\Gamma_t = \mu + t \Gamma$, depending on $t$, in the resolvent set of $\ron{\tilde{L}}$, and the integral becomes 
  \[ e^{t \ron{\tilde{L}}} = 
          \frac{-e^{-\mu t}}{2 \pi i} \, \int_{\Gamma_t} e^z \, {\cal R}\left(\frac{z}{t} - \mu\right) \, \frac{dz}{t}. \]
By the Cauchy theorem, we deform $\Gamma_t$ into some fixed $\Gamma'$, independent of $t$, 
running from infinity with $\arg(\xi)=\pi/2+\eta$ to infinity 
with $\arg(\xi) = -\pi/2 - \eta$ in the set 
  \[  \left\{ \xi \in \C, \ \ \arg (\xi) \in 
                     \left[-\frac{3\pi}{4},\frac{3\pi}{4} \right] \ \ \mbox{ and } \ \  
                     \mbox{Re}(\xi) \le \frac{\mu}{2} \right\}, \]
and the formula for the semi-group becomes 
  \[ e^{t \ron{\tilde{L}}} = 
          \frac{-e^{-\mu t}}{2 \pi i} \, \int_{\Gamma'} e^z \, 
                 {\cal R}\left(\frac{z}{t} - \mu\right) \, \frac{dz}{t}. \]
Then for any $t \ge 1$, $\Gamma'/t - \mu \subset {\cal A}_\mu$ and thus we can apply the estimate 
of Lemma~\ref{lem:resolvent} to get 
  \[ \big\| e^{t \ron{\tilde{L}}} \big\|_{L^1} = 
          \left\| \frac{-e^{-\mu t}}{2 \pi i} \, \int_{\Gamma'} e^z \, {\cal R}\left(\frac{z}{t} - \mu\right) 
              \, \frac{dz}{t} \right\|_{L^1}  
                  \le \frac{e^{-\mu t}}{2 \pi} \, \left[ a \, \int_{\Gamma'} |e^z| \, |dz| + 
                          b \, \int_{\Gamma'} |e^z| \, \frac{|dz|}{|z|} \right], \]
which concludes the proof.
\end{proof}

\Remarks 

1. This proof shows in fact that $\ron{\tilde{L}}$ is a sectorial 
operator, which implies that its evolution semi-group is analytic in $t$ 
(see~\cite[Chapter~1, Section~3]{Henry}). 
\smallskip

2. With the same method one can also prove that $\ron{L}$ is sectorial, and 
define its analytic semi-group $e^{t  \ron{L}}$ on $L^1$ which satisfies  
  \[ \forall \, t \ge 0, \quad \big\| e^{t\ron{L}} \big\|_{L^1} \le C \]
for some explicit constant $C$ depending only the collision kernel. 
More precisely if $\Pi_0$ denotes the spectral projection associated with 
the $0$ eigenvalue (for the definition of the spectral projection  
we refer to \cite[Chapter 3, Section 6, Theorem 6.17]{Kato}), then 
we have the following relation:
  \[ \forall \, t \ge 0, \quad e^{t\ron{L}} = \Pi_0 + e^{t\ron{\tilde{L}}} (\mbox{Id} - \Pi_0). \]  
\medskip


\subsection{Proof of the convergence}
\label{subsec:connect}

In this subsection we shall complete the proof of Theorem~\ref{theo:cvg}. We decompose 
the argument into several lemmas. 
The first technical lemma deals with the bilinear term $\Gamma$. 
  \begin{lemma}\label{lem:Gamma} 
  Let $B$ be a collision kernel satisfying  
  assumptions~\eqref{eq:prod}, \linebreak \eqref{eq:hyprad}, \eqref{eq:grad}. 
  Then there is an explicit constant $C_{11}>0$ 
  depending on the collision kernel such that the bilinear operator $\Gamma$ satisfies 
    \begin{equation*}
    \| \Gamma (g,g) \|_{L^1} \le C_{11} \, 
    \| g \|_{L^1} ^{3/2} \, \|g\|_{L^1(m^{-1})} ^{1/2}. 
    \end{equation*} 
  \end{lemma}
\begin{proof}[Proof of Lemma~\ref{lem:Gamma}]
Estimates in $L^1$ of the collision operator (for instance see~\cite[Section~2]{MoVi}), 
plus the obvious control 
 \[ \left( m^{-1} m' m'_* \right),\  \left( m^{-1} m m_* \right) \le 1, \]
yield 
 \[ \| \Gamma (g,g) \|_{L^1} \le C \, \|g\|_{L^1} \, \|g\|_{L^1 (\langle v \rangle^\gamma)}. \]
Then by H\"older's inequality
 \[ \|g\|_{L^1 (\langle v \rangle^\gamma)} \le C \, \|g\|_{L^1} ^{1/2} \, \|g\|_{L^1(m^{-1})} ^{1/2}. \]
\end{proof}

In a second technical lemma we state a precise form of 
the Gronwall estimate (for which we do not search for an optimal 
statement). 
  \begin{lemma}\label{lem:Gron}
  Let $y=y(t)$ be a nonnegative continuous 
  function on $\R_+$ such that 
  for some constants $a$, $b$, $\theta$, $\mu>0$, 
    \begin{equation}\label{eq:gron}
    y(t) \le a \, e^{-\mu t} y(0) + b \, 
     \left( \int_0 ^t e^{-\mu(t-s)} y(s)^{1+\theta} \, ds \right). 
    \end{equation}
  Then if $y(0)$ and $b$ are small enough, we have  
    \begin{equation*}
    y(t) \le C_{12} \, y(0) \, e^{- \mu t}.
    \end{equation*}
  for some explicit constant $C_{12}>0$. 
  \end{lemma}
\begin{proof}[Proof of Lemma~\ref{lem:Gron}]
As $y$ is continuous on $\R_+$, 
the right hand side in~\eqref{eq:gron} is continuous with respect to $t$.  
If we assume that $y(0) <1$, this remains true on a small time interval $[0,t_0]$, 
on which we have 
  \[ y(t) \le a \, e^{-\mu t} y(0) + b \, 
   \left( \int_0 ^t e^{-\mu(t-s)} y(s) \, ds \right), \]
which implies 
  \begin{equation}\label{eq:groninter}
  y(t) \le a \, y(0) \, e^{-(\mu-b) t}. 
  \end{equation}
So if we choose $y(0)$ small enough such that $a \, y(0) <1$ and 
$b$ small enough such that $\mu - b \ge 0$, we get 
for all time that $y(t) < 1$ with the bound~\eqref{eq:groninter}. 
Now to obtain the rate of decay $\mu$, we assume, by taking 
$b$ small enough, that 
 \[ \mu - b \ge \frac{\mu+\eta}{1+\theta} \]
for some $\eta >0$. We deduce that 
 \[ e^{\mu t} \, y(t) \le a \, y(0) + b \, \left( a \, y(0) \right)^{1+\theta} \,  
  \left( \int_0 ^t e^{-\eta s} \, ds \right) \le C y(0), \]
which concludes the proof. 
\end{proof}

Now we state the result of convergence to equilibrium assuming a uniform smallness estimate on 
the $L^1(m^{-1})$ norm of $g$ ({\it i.e.} the $L^1(m^{-2})$ of $f-M$).
  \begin{lemma}\label{lem:cvgpert}
  Let $B$ be a collision kernel satisfying  
  assumptions~\eqref{eq:prod}, \linebreak \eqref{eq:hyprad}, \eqref{eq:grad},  
  and $\lambda$ be the associated spectral gap. 
  Let $0<\mu \le \lambda$. Then there are some explicit constants $\var$, $C_{13}>0$ depending on the 
  collision kernel, on $\mu$ and on a lower bound on $\nu_0 - \mu$, such that if 
  $f \ge 0$ is a solution to the Boltzmann equation such that 
    \[ \forall \, t \ge 0, \ \ \ \| f_t - M \|_{L^1(m^{-2})} \le \var, \]
  then  
    \[ \forall \, t \ge 0, \ \ \ \| f_t - M \|_{L^1(m^{-1})} 
           \le C_{13} \, \| f_0 - M \|_{L^1(m^{-1})} \, e^{-\mu t}. \]
  \end{lemma}
\begin{proof}[Proof of Lemma~\ref{lem:cvgpert}]
We write a Duhamel representation of $g_t$:  
  \[ g_t = e^{t\ron{\tilde{L}}} \, g_0 + \int_0 ^t e^{(t-s)\ron{\tilde{L}}} \, \Gamma(g_s,g_s) \, ds. \]
Using Theorem~\ref{theo:sg} and Lemma~\ref{lem:Gamma} we get 
  \[ \| g_t \|_{L^1} \le C_{10} \, e^{-\mu t} \, \| g_0 \|_{L^1} 
                       + C_{10} \, C_{11} \, \var^{1/2} 
       \, \int_0 ^t e^{-\mu(t-s)} \|g_s\|_{L^1}^{3/2} \, ds. \]
Thus if $\var$ is small enough, 
we can apply Lemma~\ref{lem:Gron} with $y(t)=\|g_t\|_{L^1}$ and $\theta=1/2$ to get 
  \[ \| g_t \|_{L^1} \le C_{13} \,  \| g_0 \|_{L^1} \, \, e^{-\mu t}. \]
This concludes the proof since 
  \[ \|g_t\|_{L^1} = \|f_t -M\|_{L^1(m^{-1})}. \]
\end{proof} 

Finally we need a result on the appearance and propagation of the $L^1(m^{-1})$ norm. 
This lemma is a variant of results in~\cite{Bo:mts} and~\cite{BoGaPa}, 
and it is a particular case of more general results in~\cite{MMR1}. 
  \begin{lemma}\label{lem:mts}
  Let $B$ be a collision kernel satisfying  
  assumptions~\eqref{eq:prod}, \eqref{eq:hyprad}, \eqref{eq:grad}, \eqref{eq:hypmts}.
  Let $f_0$ be a nonnegative initial datum in $L^1(\langle v \rangle^2) \cap L^2$. 
  Then the corresponding solution $f=f(t,v)$ of the Boltzmann equation~\eqref{eq:base}  
  in $L^1(\langle v \rangle^2)$ satisfies: for any $0<s<\gamma/2$ and $\tau >0$, there 
  are explicit constants $a, C >0$, depending on the collision kernel, $s$,  
  $\tau$, and the mass and energy and $L^2$ norm of $f_0$, such that 
     \begin{equation}\label{eq:mts}
     \forall \, t \ge \tau, \ \ \ \int_{\R^N} f(t,v) \, \exp\left[a |v|^s \right] \, dv  \le C. 
     \end{equation}
  In the important case of hard spheres~\eqref{eq:hs}, the assumption 
  ``$f_0 \in L^1(\langle v \rangle^2) \cap L^2$'' 
  can be relaxed into just ``$f_0 \in L^1 (\langle v \rangle^2)$'', 
  and the same result~\eqref{eq:mts} holds 
  with explicit constant $a,C >0$ depending only on the collision 
  kernel, $s$, $\tau$, and the mass and energy of $f_0$. 
  \end{lemma}
\begin{proof}[Proof of Lemma~\ref{lem:mts}]
Note that the assumption $f_0 \in L^2$ implies in particular 
that $f_0$ has finite entropy, {\it i.e.} 
  \[ \int_{\R^N} f_0(v) \log f_0(v) \, dv \le H_0 < +\infty\]
which ensures by the $H$ theorem that 
  \begin{equation}\label{nonconc}
  \forall \, t \ge 0, \ \ \ \int_{\R^N} f(t,v) \log f(t,v) \, dv \le H_0.
  \end{equation}
We assume, up to a normalization, that $f$ satisfies 
  \begin{equation}\label{centre}
  \forall \, t \ge 0, \ \ \ \int_{\R^N} f(t,v) \, v \, dv = 0. 
  \end{equation}

Let us fix $0<s<\gamma/2$. We define for any $p \in \R_+$
  \[ m_p(t) := \int_{\R^N} f(t,v) \, |v|^{sp} \, dv. \]
The evolution equation on the distribution $f$ yields 
  \begin{equation} \label{mp1}
  \frac{d m_p}{dt} = \int_{\R^N} Q(f,f) \, |v|^{sp} \, dv
  =  \int_{\R^N \times \R^N} f \, f_* \, \Phi(|v-v_*|) \,  
  K_p (v,v_*) \, dv \, dv_*,
  \end{equation}
where
  \begin{equation} \label{mp2}
  K_p(v,v_*) := \frac12 \, \int_{S^{N-1}} \big(|v'|^{sp} + |v'_*|^{sp} - |v|^{sp} - 
  |v_*|^{sp}\big) \, b(\cos \theta) \, d\sigma.
  \end{equation}
From~\cite[Lemma 1, Corollary 3]{BoGaPa}, we have
  \begin{equation} \label{mp3}
  K_p(v,v_*) \le \alpha_p \, \big(|v|^2 + |v_*|^2 \big)^{sp/2} - |v|^{sp} - |v_*|^{sp}
  \end{equation}
where $(\alpha_p)_{p \in \N/2}$ is a strictly decreasing sequence 
such that
  \begin{equation} \label{mp4}
  0< \alpha_p < \min \left\{ 1, \frac{C}{sp/2+1} \right\}.
  \end{equation}
for some constant $C$ depending on $C_b$, defined in~\eqref{eq:grad}.  
Notice that the assumptions~\cite[(2.11)-(2.12)-(2.13)]{BoGaPa} 
are satisfied under our assumptions \eqref{eq:prod}, 
\eqref{eq:hyprad}, \eqref{eq:grad}, \eqref{eq:hypmts} on the collision kernel 
(see \cite[Remark~3]{BoGaPa} for 
some possible ways of relaxing the assumption \eqref{eq:hypmts} in the elastic case).     

Then we use the classical estimate 
  \[ \int_{\R^N} \Phi(|v-v_*|) \, f(t,v_*) \, dv_* 
     = C_\Phi \, \int_{\R^N} |v-v_*|^\gamma \, f(t,v_*) \, dv_* \ge K \, |v|^\gamma \]
for some constant $K$ depending on $C_\Phi$ (defined in~\eqref{eq:hyprad}) and the 
mass, energy and entropy $H_0$ of $f_0$ (or only the mass of $f_0$ in the case $\gamma=1$). \
This estimate is obtained by a classical non-concentration argument when $\gamma \in (0,1)$ 
(using the bound~\eqref{nonconc}), as can be found in~\cite{Ar72} for instance,  
or just by convexity when $\gamma=1$, using the assumption~\eqref{centre} that 
the distribution has zero mean (see for instance~\cite[Lemma~2.3]{MMR1}). 
We combine this with~\cite[Lemma~2 and Lemma~3]{BoGaPa} to get 
  \begin{equation} \label{mp5}
  \int_{\R^N} Q(f,f) \, |v|^{sp} \, dv
  \le \alpha_p \, Q_p - K \, (1-\alpha_p) \, m_{p+\gamma/s}
  \end{equation}
with
  \[ Q_p := \int_{\R^N \times \R^N} f f_* \left[ \big(|v|^2 + |v_*|^2 \big)^{sp/2} 
                               - |v|^{sp} - |v_*|^{sp} \right] \, \Phi(|v-v_*|) \, dv \, dv_*, \]
and~\cite[Lemma~2 and Lemma~3]{BoGaPa} shows that, for $ps/2 > 1$, we have 
  \begin{equation}\label{mp5bis}
  Q_p \le  S _p = C_\Phi \, \sum_{k=1}^{k_p} 
  \left(\begin{matrix} sp/2 \\ k  \end{matrix} \right) \big(m_{(2k+\gamma)/s} \, m_{p-2k/s} + m_{2k/s} 
  \, m_{p-2k/s+\gamma/s}\big), 
  \end{equation}
where $k_p := [sp/4+1/2]$ is the integer part of $(sp/4+1/2)$ and 
$\left(\begin{smallmatrix} sp/2 \\ k  \end{smallmatrix} \right)$ denotes the 
generalized binomial coefficient. Gathering~\eqref{mp1}, \eqref{mp5} and~\eqref{mp5bis}, we get
  \begin{equation} \label{mp6}
  \frac{dm_p}{dt} \le \alpha_p \, S _p - K \, (1-\alpha_p) \,  m_{p+\gamma/s}. 
  \end{equation}
By H\"older's inequality, we have 
  \[ m_{p+\gamma/s} \ge \beta \, m_p^{1+\frac{\gamma}{sp}}  \]
for some constant $\beta >0$ depending on the mass of the distribution.  
By~\cite[Lemma~4]{BoGaPa}, there exists $A > 0$ such that
  \[ S _p \le A \, \Gamma (p+1+\gamma/s) \, Z_p \]
with
  \[ Z_p := \max_{k=1,..,k_p} \big\{ z_{(2k+\gamma)/s} \, z_{p-2k/s}, \,  z_{2k/s} \, 
     z_{p-2k/s+\gamma/s} \big\}, \ \ \mbox{ and } \ \ z_p := \frac{m_p}{\Gamma(p+1/2)}. \]
Thus we may rewrite~\eqref{mp6} as
  \begin{equation} \label{mp7}
  \frac{d z_p}{dt} \le 
  A \, \alpha_p \frac{\Gamma(p+1 + \gamma/s)}{\Gamma(p+1/2)} \, Z_p
  - K' \, (1-\alpha_p) \, \Gamma(p+1/2)^{\frac{\gamma}{sp}} \, z_p^{1+\frac{\gamma}{sp}}
  \end{equation}
with $K' = \beta \, K$. On the one hand, from the definition of the sequence  
$(\alpha_p)_{p \ge 0}$, there exists $A'$ such that
  \begin{equation} \label{mp8}
  A \, \alpha_p \, \frac{\Gamma(p+1+\gamma/s)}{\Gamma(p+1/2)}  \le A' \, p^{\gamma/s-1/2}.
  \end{equation}
On the other hand, thanks to Stirling's formula 
  \[ n! \begin{array}{c} \vspace{-0.1cm} \\ \sim \vspace{-0.25cm} \\ \scriptstyle n \to +\infty \end{array}  
       n^n \, e^{-n} \, \sqrt{2\pi n}, \]
there is $A'' > 0$ such that
  \begin{equation} \label{mp9}
  (1-\alpha_p) \, \Gamma(p+1/2)^{\frac{\gamma}{sp}} \ge A'' \, p^{\gamma/s}. 
  \end{equation}
Gathering~\eqref{mp7}, \eqref{mp8} and~\eqref{mp9}, we deduce
  \begin{equation}\label{mp10}
  \frac{d z_p}{dt} \le A' \, p^{\gamma/s -1/2} \, Z_p 
     - A''  \, K' \, p^{\gamma/s} \, z_p^{1+\frac{\gamma}{sp}}. 
  \end{equation}
Note that in this differential inequality, for $p$ big enough, 
$Z_p$ depends only on terms $z_q$ for $q \le p-1$, (but not 
necessarily with $q$ an integer), which allows 
to get bounds on the moments by induction. 

At this stage, we prove by induction on $p \ge p_0$ integer ($p_0 \ge 1$) the following property
  \[ \forall \, t \ge t_p, \ \forall \, q \in [p_0,p]\ \ \ z_q \le x^q  \]
for some $x \in (1, +\infty)$ large enough and some increasing sequence of times 
$(t_p)_{p \ge p_0}$ with $t_{p_0} >0$ fixed as small as wanted. 
The goal is to prove this induction for a convergent 
sequence of times $(t_p)_{p \ge p_0}$. 
The initialization for $p = p_0$, with $p_0$ as big as wanted and $t_0$ 
as small as wanted, is straightforward by the classical 
theorems about the immediate appearance and uniform propagation of algebraic moments 
(see \cite{Wenn:momt:97} for instance), and taking $x$ big enough. 
Now as $s <\gamma/2$, if $p_0$ is large enough, we have for $p \ge p_0$
  \[ p^{\gamma/s} \ge 2 \, \frac{A'}{A'' \, K'} \, p^{\gamma/s-1/2},  
     \qquad  p^{\gamma/s} \ge p^{2 + \var} \ \mbox{ and } \ p \ge \left( \frac{A'}{A'' K'} \right)^2 \]
for some $\var >0$. So let us assume the induction property satisfied for all 
steps $p_0 \le q \le p-1$, and let us consider $z_p$. 
Assume that $z_p(t_{p-1}) \le x^p$.  
Then from~\eqref{mp10}, for any $t$ such that $z_p(t) \le x^p$ we have 
  \begin{eqnarray*}
  \frac{d z_p}{dt} &\le& A' \, p^{\gamma/s -1/2} \, Z_p 
     - A''  \, K' \, p^{\gamma/s} \, z_p^{1+\frac{\gamma}{sp}} \\
  &\le& A' \, p^{\gamma/s -1/2} \, \left[ (x^p)^{1+\frac{\gamma}{sp}} 
         - \frac{A'' K'}{A'} \, p^{1/2} \, z_p^{1+\frac{\gamma}{sp}} \right] \\
  &\le& A' \, p^{\gamma/s -1/2} \, \left[ (x^p)^{1+\frac{\gamma}{sp}} - z_p^{1+\frac{\gamma}{sp}} \right].
  \end{eqnarray*}
We deduce by maximum principle 
that this bound is propagated uniformly for all times $t \ge t_{p-1}$. 
If on the other hand $z_p(t_{p-1}) \in (x^p,+\infty]$, as long as $z_p(t) > x^p$ we have
  \begin{eqnarray*}
  \frac{d z_p}{dt} &\le& A' \, p^{\gamma/s -1/2} \, Z_p 
     - A''  \, K' \, p^{\gamma/s} \, z_p^{1+\frac{\gamma}{sp}}\\
     &\le&  A' \, p^{\gamma/s -1/2} \, (x^p)^{1+\frac{\gamma}{sp}} 
     - A''  \, K' \, p^{\gamma/s} \, z_p^{1+\frac{\gamma}{sp}} \\
     &\le& \left[ A' \, p^{\gamma/s -1/2} - A''  \, K' \, p^{\gamma/s} \right] \, z_p^{1+\frac{\gamma}{sp}} \\
     &\le& - \frac{A''  \, K'}{2} \, p^{\gamma/s} \, z_p^{1+\frac{\gamma}{sp}} 
     \le -C \, p^{2 + \var} \,  z_p^{1+\theta}. 
  \end{eqnarray*}
with $C = (A''  \, K')/2$ and $\theta = \gamma/(sp)$. Then classical arguments of 
comparison to a differential equation show that $z_p$ is finite 
for any $t > t_{p-1}$, and satisfies the following bound independent of the initial datum:  
  \[ z_p ( t_{p-1} + t ) \le \left[ \frac{1}{C \, \theta \, p^{2+\var} \, t} \right]^{\frac{1}{\theta}}
                        \le \left[ \frac{s}{C \, \gamma \, p^{1+\var} \, t} \right]^{\frac{sp}{\gamma}}. \]
Thus if we set 
  \[ \Delta_p :=  \frac{s}{C \, \gamma \, p^{1+\var} \, x^{\gamma/s}}, \]
we have 
  \[ \forall \, t \ge \Delta_p, \quad z_p ( t_{p-1} + t ) \le x^p. \]
By H\"older interpolation with the bounds on $z_{p-1}$ 
provided by the induction assumption, we deduce immediately 
  \[ \forall \, t \ge \Delta_p, \ \forall \, q \in [p,p+1], \quad z_q ( t_{p-1} + t ) \le x^q, \]
which proves the step $p$ of the induction with $t_p = t_{p-1} + \Delta_p$. 
So we have proved that if we set  
  \[ \tau = \lim_{p \to +\infty} t_p = t_{p_0} + \sum_{p \ge p_0+1} \Delta_p =  
       t_0 + \frac{s}{C \, \gamma \, x^{\gamma/s}} 
                   \, \left( \sum_{p=p_0+1} ^{+\infty} \frac{1}{p^{1+\var}} \right) < +\infty, \]
we have 
  \[ \forall \, t \ge \tau, \ \forall \, p \ge p_0, \quad z_p (t) \le x^p. \]
Moreover $\tau$ can be taken as small as wanted by taking $t_0$ small enough and 
$x$ large enough. 

So we conclude that there are explicit constants $R$ and $\tau$ such that
  \[ \forall \, t \ge \tau, \ \forall \, p \ge 0,\ \ \  z_p(t) \le R^{-p}. \]
We deduce explicit bounds
  \[ \sup_{t \ge \tau} \int_{\R^N} f(t,v) \, \exp \left[ a |v|^s \right] \, dv \le C <+\infty \]
for any $a < R$ since 
  \begin{multline*}
  \int_{\R^N} f(t,v) \, \exp \left[ a |v|^s \right] \, dv = \sum_{p=0}^\infty 
  \int_{\R^N}f(t,v) \, a^p \, \frac{|v|^{sp}}{p!} \, dv =
  \sum_{p=0}^\infty a^p \, \frac{m_p}{p!} \\
  = \sum_{p=0}^\infty a^p \, \frac{z_p \, \Gamma(p+1/2)}{p!} \le
  \sum_{p=0}^\infty \left( \frac{a}{R} \right)^p \frac{\Gamma(p+1/2)}{p!}  
  \le C \, \sum_{p=0}^\infty \left( \frac{a}{R} \right)^p \, p^{1/2}  < +\infty.
  \end{multline*}
\end{proof}

Now we state a result of convergence to equilibrium for non smooth 
solutions deduced from~\cite[Theorems~6.2 and~7.2]{MoVi} combined 
with the previous lemma. 
  \begin{lemma}\label{lem:mv}
  Let $B$ be a collision kernel satisfying  
  assumptions~\eqref{eq:prod}, \eqref{eq:hyprad}, \eqref{eq:grad}, \eqref{eq:hypmts}, 
  \eqref{eq:beep}. Let $f_0$ be a nonnegative initial datum in $L^1(\langle v \rangle^2) \cap L^2$. 
  Then the corresponding solution $f \ge 0$ of the Boltzmann equation~\eqref{eq:base}  
  in $L^1(\langle v \rangle^2)$ satisfies: for any $\tau>0$, $0 < s < \gamma/2$, there 
  are explicit constants $C_{14} >0$ and $a>0$ depending on the collision kernel,  
  $\tau$, $s$ and the mass, energy and $L^2$ norm of $f_0$, such that
    \begin{equation}\label{eq:mv}
    \forall \, t \ge \tau, \ \ \ \|f_t -M\|_{L^1 (m^{-1})} 
                \le  C_{14} \, t^{-1} 
    \end{equation}
  with $m(v) = \exp \left[-a|v|^s \right]$.  
  In the important case of hard spheres~\eqref{eq:hs}, the assumption ``$f_0 \in L^1(\langle v \rangle^2) \cap L^2$'' 
  can be relaxed into just ``$f_0 \in L^1(\langle v \rangle^2)$'', and the same result~\eqref{eq:mv} holds 
  with explicit constant $a,C_{14} >0$ depending only on the collision 
  kernel, $s$, $\tau$, and the mass and energy of $f_0$. 
  \end{lemma}
 
\begin{proof}[Proof of Lemma~\ref{lem:mv}]
It is straightforward that the assumptions~\eqref{eq:prod}, \eqref{eq:hyprad}, \eqref{eq:grad}, 
\eqref{eq:beep} implies the assumptions~\cite[equations~(1.2) to~(1.7)]{MoVi}. 
Hence by~\cite[Theorem~6.2]{MoVi} we deduce that for an initial datum in $L^1 (\langle v \rangle^2) \cap L^2$, 
the solution satisfies for any $\alpha >0$ 
  \begin{equation}\label{eq:mvprec}
  \forall \, t \ge 0, \ \ \  \|f_t - M\|_{L^1} 
                \le C_\alpha \, t^{-\alpha} 
  \end{equation}
for some explicit constant $C_\alpha$ depending on $\alpha$, the collision kernel 
and the mass, energy and $L^2$ norm of the initial datum.  
We apply this result with $\alpha=2$ and we interpolate by H\"older's inequality 
with the norm $L^1 (\exp\left[2a |v|^s\right])$ for $0<s<\gamma/2$ which is bounded 
uniformly for $t \ge \tau>0$ by Lemma~\ref{lem:mts}, to deduce that 
  \[  \forall \, t \ge \tau, \ \ \  \|f_t - M\|_{L^1 (m^{-1})} 
                \le C_{14} \, t^{-1} \]
with $m(v)= \exp\left[ - a |v|^s \right]$. 

In the case of hard spheres we use~\cite[Theorem~7.2]{MoVi} instead of~\cite[Theorem~6.2]{MoVi}, 
which yields the same result~\eqref{eq:mvprec}, but under the sole assumption on the initial 
datum that $f_0 \in L^1 _2$. 
\end{proof} 

Now we can conclude the proof of Theorem~\ref{theo:cvg}:
\begin{proof}[Proof of Theorem~\ref{theo:cvg}]
Using Lemmas~\ref{lem:mts} and~\ref{lem:mv}, 
we pick $t_0>0$ and $m(v) = \exp \left[ -a|v|^s \right]$ with $0< s < \gamma/2$ such that 
  \[  \forall \, t \ge t_0, \ \ \ \|f_t -M\|_{L^1(m^{-2})} \le \var \]
where $\var$ is chosen as in Lemma~\ref{lem:cvgpert}. Then for $0<\mu\le \lambda$, we apply 
Lemma~\ref{lem:cvgpert} starting from $t=t_0$: 
  \[ \forall \, t \ge t_0, \ \ \ \| f_t - M \|_{L^1(m^{-1})} 
              \le C_{13} \, \| f_{t_0} - M \|_{L^1(m^{-1})} \, e^{-\mu t}
              \le C \, e^{-\mu t}. \]
This concludes the proof. 
\end{proof}

\subsection{A remark on the asymptotic behavior of the solution}
\label{subsec:profil}

Theorem~\ref{theo:cvg} thus yields the asymptotic expansion 
 \[ f = M + m g \] 
with $g$ going to $0$ in $L^1$ with rate $C \, e^{-\lambda t}$. 
If we denote by $\Pi_1$ the spectral projection associated with the 
$-\lambda$ eigenvalue (for the definition of the spectral projection  
see \cite[Chapter 3, Section 6, Theorem 6.17]{Kato}) and 
$\Pi_1 ^\bot = \mbox{Id} - \Pi_1$, then the evolution equation on 
$\Pi_1 ^\bot(g)$ writes (using the fact that $\ron{L}$ commutes with $\Pi_1$)  
 \[ \partial_t \Pi_1 ^\bot(g) = \ron{L} ( \Pi_1 ^\bot(g)) + \Pi_1 ^\bot(\Gamma(g,g)). \]
By exactly the same analysis as above, one could prove that the semi-group of 
$\ron{L}$ restricted to $\Pi_1 ^\bot(L^1)$ decays with rate 
$C \, e^{-\lambda_2 t}$ where $\lambda_2>\lambda$ is the modulus of the 
second non-zero eigenvalue. Thus by the Duhamel formula one gets 
 \[ \| \Pi_1 ^\bot(g_t) \|_{L^1}  \le C \, e^{-\lambda_2 \, t} \, \| \Pi_1 ^\bot(g_0) \|_{L^1} 
                 + C \, \int_0 ^t e^{-\lambda_2 (t-s)} \, \|\Pi_1 ^\bot(\Gamma(g_s,g_s))\|_{L^1} \, ds. \]
Then using that 
 \[  \|\Pi_1 ^\bot(\Gamma(g_s,g_s))\|_{L^1} \le  C \, \| \Gamma(g_s,g_s)\|_{L^1} 
               \le C \, \|g_s \|_{L^1} ^{3/2} 
                \le C \, e^{-(3/2) \lambda t} \]
we deduce that 
 \[ \| \Pi_1 ^\bot(g_t) \|_{L^1} \le C \, e^{-\bar{\lambda} t} \]
with $\lambda < \bar{\lambda} < \min \{ \lambda_2, (3/2) \lambda \}$. Hence, 
setting $\varphi_1 = m \, \Pi_1(g_t)$ and $R= m \, \Pi^\bot(g_t)$, we obtain the 
asymptotic expansion
 \[ f = M + \varphi_1 + R \]
with $\varphi_1=\varphi_1(t,v)$ going to $0$ in $L^1$ with rate $e^{-\lambda t}$ and 
$R=R(t,v)$ going to $0$ in $L^1$ with rate $e^{-(\lambda + \var ) t}$ for some 
$\var > 0$. Thus $\varphi_1$ is asymptotically the dominant term of $f-M$, and as 
$m^{-1} \varphi_1$ belongs to the eigenspace of $\ron{L}$ associated with $\lambda$, 
we know by the study of the decay of the eigenvectors that 
 \[ \forall \, t \ge 0, \hspace{0.3cm} \varphi_1  \in L^2 (M^{-1}). \]
Moreover $\varphi_1$ is the projection of the solution on the eigenspace of the 
first non-zero eigenvalue. It can be seen as the first order correction to 
the equilibrium regime, and the asymptotic profil of this first order 
correction is given by the eigenvectors associated to the $-\lambda$ eigenvalue.  
\bigskip

\noindent {\bf{Acknowledgments.}} 
The idea of searching for some decay property on the 
eigenvectors imposed by the eigenvalue equation in order to 
prove that the eigenvectors belong to a smaller space of 
linearization originated from fruitful discussions 
with Thierry Gallay, under the impulsion  
of C\'edric Villani. Both are gratefully acknowledged. 
We also thank St\'ephane Mischler for useful remarks on the 
preprint version of this work. 
Support by the European network HYKE, funded by the EC as
contract HPRN-CT-2002-00282, is acknowledged.

\begin{flushright} \signcm \end{flushright}
\end{document}